%% file: IRO7.tex
\documentclass[12pt]{article}

\usepackage{multicol}

\usepackage{amssymb}
\usepackage[leqno]{amsmath}
\usepackage{amsthm}

\usepackage{graphicx,color}

\usepackage{eso-pic}
\definecolor{lightgray}{gray}{.85}
\usepackage{subfigure}

\usepackage{xr-hyper}
\usepackage[bookmarks,pdfnewwindow]{hyperref}

\hypersetup{
colorlinks=true,
linkcolor=black,
citecolor=black,
pdfauthor={Victor Ivrii},
pdftitle={Sharp Spectral Asymptotics for 2-dimensional Schr\"odinger operator with a strong but degenerating magnetic field. II},
pdfsubject={Sharp Spectral Asymptotics},
pdfkeywords={Microlocal Analysis, Magnetic Schr\"odinger Operator},
baseurl={http://www.math.toronto.edu/ivrii/Research/preprints/},
}

\newcommand\W{{\rm W}}
\newcommand\MW{{\rm {MW}}}
\newcommand\Def{{\overset {\rm {def}}{\ =\ }}}

\newcommand\out{{\rm {out}}}
\newcommand\inn{{\rm {inn}}}

\newcommand\const{{\rm {const}}}

\newcommand\eff{{\rm {eff}}}
\newcommand\corr{{\rm {corr}}}

\newcommand\Tr{\operatorname{Tr}}

\newcommand\bR{{\mathbb R}}

\newcommand\bH{{\mathbb H}}

\newcommand\bZ{{\mathbb Z}}
\newcommand\cA{{\mathcal A}}
\newcommand\cB{{\mathcal B}}

\newcommand\cE{{\mathcal E}}

\newcommand\cR{{\mathcal R}}

\newcommand\cU{{\mathcal U}}

\newcommand\cX{{\mathcal X}}
\newcommand\cY{{\mathcal Y}}
\newcommand\cZ{{\mathcal Z}}

\newtheorem{theorem}{Theorem}[section]

\newtheorem{corollary}[theorem]{Corollary}
\newtheorem{proposition}[theorem]{Proposition}

\theoremstyle{definition}

\newtheorem{remark}[theorem]{Remark}

\newcommand{\sect}[1]{\setcounter{equation}{0}\section{#1}}

\renewcommand{\theequation}{\thesection.\arabic{equation}}

\newenvironment{claim}[1][(\theequation)]{\refstepcounter{equation}\smallskip
\begin{trivlist}
\item[{\hskip\labelsep#1}]}{\smallskip\end{trivlist}}

\newcounter{note}

\newenvironment{claim*}[1]{\smallskip
\begin{trivlist}
\item[{\hskip\labelsep#1}]}{\smallskip\end{trivlist}}

\newenvironment{pdeq}[1][]{\refstepcounter{equation}}{}

\begin{document}

\title{%
Sharp Spectral Asymptotics for 2-dimensional Schr\"odinger operator with a strong
but degenerating magnetic field. II}

\author{%
Victor Ivrii
\footnote{Work was partially supported by NSERC grant OGP0138277.}
}
\maketitle

{\abstract
I consider the same operator as in part I \cite{IRO6} assuming however that
$\mu \ge Ch^{-1}$ and $V$ is replaced by $(2l+1)\mu h F+W$ with $l\in \bZ^+$. Under some non-degeneracy conditions I recover remainder estimates up to $O \bigl( \mu^{-{\frac 1 \nu}}h^{-1} +1\bigr)$ but now case $\mu \ge Ch^{-\nu}$ is no more forbidden and the principal part is of magnitude $\mu h^{-1}$. 

\endabstract}


\setcounter{section}{5}
\sect{Modified $V$. I. $\mu \le \epsilon \lowercase {h}^{-\nu}$}
\subsection{Introduction}

This paper is a continuation of \cite{IRO6} which is considered as Part I. I consider spectral asymptotics of the magnetic Schr\"odinger operator
\begin{equation}
A= {\frac 1 2}\Bigl(\sum_{j,k}P_jg^{jk}(x)P_k -V\Bigr),\qquad P_j=D_j-\mu V_j,\quad V= (2l+1)\mu h F+W
\label{6-1}
\end{equation}
where $g^{jk}$, $V_j$, $W$ are smooth real-valued functions of $x\in \bR^2$, $l\in \bZ^+$ (i.e. $l=0,1,\dots$) and
$(g^{jk})$ is positive-definite matrix, $0<h\ll 1$ is a Planck parameter and
$\mu \gg1$ is a coupling parameter. I assume that $A$ is a self-adjoint operator and all the conditions are satisfied in the ball $B(0,1)$,
$F=F_{12}g^{-{\frac 1 2}}$, $F_{12}=\partial_{x_1}V_2-\partial_{x_2}V_1$,
$g=\det (g^{jk})^{-1}$. 

Further, exactly as in \cite{IRO6}, I assume that
\begin{equation}
F \asymp |x_1|^{\nu -1},\qquad \nu\in \bZ^+,\ \nu \ge 2
\label{6-2}
\end{equation}
and thus with no loss of the generality I can assume that
\begin{equation}
V_1=0,\qquad V_2\asymp x_1^\nu.
\label{6-3}
\end{equation}
Furthermore,  I assume that \underline{either}
\begin{pdeq}\label{6-4}\end{pdeq}
\begin{align}
&\pm W\ge \epsilon_0, \qquad &&\text{as\ } l\ge 0\tag*{$(\theequation)_\pm$}
\label{6-4-pm}\\
\intertext{and as $l=0$ only sign ``$+$'' is interesting \underline{or}}
&|\partial_{x_2} W/f |\ge \epsilon_0, \qquad f \Def F x_1^{1-\nu}.
\label{6-5}
\end{align}
Also as in \cite{IRO6}, I am interested in the asymptotics of
$\int e(x,x,0)\psi (x)\,dx$ where $e(x,y,\tau)$ is the Schwartz kernel of the spectral projector $E(\tau)$ of $A$ and $\psi\in C_0^\infty (B(0,{\frac 1 2}))$
is a cut-off function and I expect the main part of it to be
$\int \cE^\MW (x,0)\psi (x)\,dx$ where $\cE^\MW$ is defined by (0.8)\,\footnote{\label{foot-1} References by default are to \cite{IRO6}.}
which is of magnitude $\mu h^{-1}$. \emph{I am assuming without mention that $\psi$ is supported in the small but fixed vicinity of $\{x_1=0\}$.}

In the sharp contrast to the analysis of Part I the case $\mu \ge C h^{-\nu}$ is not ``forbidden'' anymore as well as zone
$\cZ''=\{|x_1|\ge {\bar\gamma}_1\Def C(\mu h)^{-1/(\nu -1)}\}$.
On the contrary, as $\mu h\ge C$ this zone becomes the main contributor to the principal part of asymptotics which now is of magnitude $\mu h^{-1}$ instead of
$h^{-2} (\mu h)^{-1/(\nu -1)}$ as it was in \cite{IRO6}.
Actually I will time to time slightly change the definition of ${\bar\gamma}_1$, replacing it by
${\bar\gamma}_1=\epsilon(\mu h)^{-1/(\nu -1)}$ and back and changing respectively definition of zones.

Section 6 is devoted to the case  of $\mu \le \epsilon h^{-\nu}$. Analysis in zone  $\cZ'\Def \{|x_1|\le 2{\bar\gamma}_1\}$ remains basically the same and the main attention is paid here to the formally forbidden zone $\cZ''$. The main results here are theorems \ref{thm-6-10},  \ref{thm-6-11} and \ref{thm-6-17}.

As $\mu \ge \epsilon h^{-\nu}$ this separation to zones is no more reasonable and will be modified. In section 7 I analyze the case of 
$\epsilon h^{-\nu}\le \mu \le  Ch^{-\nu}$.  The main results here are theorems \ref{thm-7-3} and  \ref{thm-7-4}.

Further, in section 8 analyze the case of  $\mu \ge  Ch^{-\nu}$. The main results here are theorems \ref{thm-8-9},  \ref{thm-8-10},  \ref{thm-8-11} and  \ref{thm-8-12}.

Finally, appendix A is devoted to asymptotics of some one-dimensional Schr\"odiunger operators associated with (\ref{6-1}).

\subsection{Simple Rescaling}

As in \cite{IRO6} the simple rescaling arguments help us to get the easy but not sharp results. 

\medskip\noindent 
{\bf \thesubsection.1 }
In this and the next subsubsection I assume that $\mu \le Ch^{-\nu}$.Rescaling arguments in the zone
$\cZ'$ work exactly in the same manner as in \cite{IRO6} leading to the asymptotics of
$\int e(x,x,0)\psi'(x)\,dx$ with the principal part
$\int \cE^\MW (x,0)\psi' (x)\,dx$ and the remainder estimate $O(h^{-1})$ where $\psi'(x)$ and $\psi''(x)$ are cut-off functions supported in zones $\cZ'$ and $\cZ''$ (defined as above) respectively; one can take
$\psi'(x)=\psi (x) \psi_0(x_1/{\bar\gamma}_1)$, $\psi''=\psi-\psi'$ where $\psi_0\in C_0^\infty $ is supported in $(-1,1)$ and equals 1 in $[-{\frac 1 2}, {\frac 1 2}]$.

However the contribution of the previously forbidden zone $\cZ''$ to the remainder estimate is
\begin{equation*}
O\Bigl(\int_{\{{\bar\gamma}_1\le \gamma \le 1\}} \gamma^{-2}\, d\gamma\Bigr)= O\bigl({\bar\gamma}_1^{-1}\bigr)
\end{equation*}
which is $O\bigl(h^{-1}\bigr)$ due to assumption $\mu \le Ch^{-\nu}$ and the contribution of $\cZ''$ to the principal part is
\begin{equation}
\int \cE^\MW (x)\psi ''(x)\,dx=
{\frac 1 {4\pi}}\mu h^{-1}l_\pm \int\psi '' |F|\sqrt g\, dx,\qquad
l_\pm \Def l+{\frac 1 2}(-1\pm 1)
\label{6-6}
\end{equation}
under condition \ref{6-4-pm}.

Under condition (\ref{6-5}) the above arguments remain true for the contribution of the subzone $\cZ''\cap \{|W|\ge C\gamma\}$; for the contribution of the zone
$\cZ''\cap \{|W|\le C\gamma\}$ one needs to take in account correction term\footnote{\label{foot-2}See section 6 of \cite{Ivr1}.}
$\sum_m\kappa_m \mu_\eff h_\eff ^{1+2m}$ for the case $\mu_\eff h_\eff\ge 1$, $h_\eff\le 1$ where in the rescaling and division arguments
$\mu_\eff=\mu \gamma^{\nu-{\frac 1 2}}$,
$h_\eff =h\gamma ^{-\frac 3 2}$ and the number of balls is $O(1)$ for each $\gamma$. Then the total contribution of this correction terms is
$O\bigl(\mu h \bigr)$ as $\nu \ge 3$ and $O\bigl(\mu h |\log h|\bigr )$ as $\nu=2$.

\medskip\noindent
{\bf \thesubsection.2 }
Replacing $\psi$ by $x_1\psi$ in the above arguments one gains factor $\gamma$ in each integrand; then the total contribution of the zone $\cZ'$ to the remainder estimate becomes
\begin{equation*}
O\Bigl(\int \mu^{-1}h^{-1}\gamma^{1-\nu}\times \gamma \times \gamma^{-2}\,d\gamma \Bigr) = O\bigl( \mu^{-1/\nu}h^{-1}\bigr)
\end{equation*}
which is exactly what I want. On the other hand, the contribution of zone $\cZ''$ to the remainder estimate becomes
$O(\gamma^{-1}\,d\gamma)=O(|\log h|)$ which is what we want as 
$\mu \le C(h|\log h|)^{-\nu}$ only. To fix it under condition \ref{6-4-pm} one can notice that zone $\cZ''$ is  the spectral gap and therefore the contribution of the individual ball to the remainder estimate is
$O(\gamma h_\eff^s)$ with $h_\eff= h/\gamma$ rather than $O(1)$ and therefore the total contribution of zone $\cZ''$ to the remainder estimate is $O(1)$.

As before, under condition (\ref{6-5}) these arguments are applicable in the subzone $\cZ''\cap \{|W|\ge C\gamma\}$ with $h_\eff= h/(\gamma |\W|^{1/2})$
as long as $h_\eff\le 1$. This leads to $O(1)$ estimate of the contribution of the subzone $\cZ'' \cap\{ |W|\ge C\gamma, |W|^{1/2}\gamma \ge h\}$ to the remainder. One can see easily that the integral of $\gamma^{-1}$ taken over subzones $\cZ'' \cap\{ |W|\ge C\gamma, |W|^{1/2}\gamma \le h\}$
and $\cZ'' \cap\{ |W|\le C\gamma\}$ is $O(1)$ as well. Thus rescaling arguments
provide remainder estimate $O\bigl(\mu^{-{1/\nu}}+1\bigr)$ if $\psi$ contains an extra factor $x_1$ and under condition (\ref{6-5}) correction terms are taken into account.

Therefore

\begin{claim}\label{6-7}
As $\mu \le Ch^{-\nu}$ in what follows one can assume without any loss of the generality that $\psi (x)= \psi_1(x_1)\psi_2(x_2)$.
\end{claim}

\medskip\noindent
{\bf \thesubsection.3 }
As $\mu \ge Ch^\nu$ arguments of subsubsection \thesubsection.1 work as $\{|x_1|\ge C h\}$ providing $O(h^{-1})$ contribution of this zone to the remainder estimate while the contribution of zone $\{|x_1|\le C h\}$ will be $O(\mu h^{\nu -1})$. The main part of the asymptotics will be the same as above.

Moreover, arguments of subsubsection \thesubsection.1 work as $\{|x_1|\ge C h\}$ providing $O(1)$ contribution of this zone to the remainder estimate as $\psi$ is replaced by $x_1\psi$
while the contribution of zone $\{|x_1|\le C h\}$ will be $O(\mu h^\nu )$.

In the next section I will improve these latter results.

\subsection{Estimates. I}
In section 2 and subsections 4.1--4.4 of \cite{IRO6} various properties of operator $A$ were proven in the outer and inner zones
$\cZ_\out =\{{\bar \gamma}\le |x_1|\le 2{\bar\gamma}_1\}$ and
$\cZ_\inn =\{|x_1|\le 2{\bar\gamma}\}$ with
${\bar\gamma}\Def C \mu^{-1/\nu}$ as long as ${\bar\gamma}\le {\bar\gamma}_1$ i.e. $\mu \le \epsilon h^{-\nu}$.
These properties were proven first in section 2 under assumption
\begin{align}
&C\le \mu \le \epsilon (h|\log h|)^{-\nu}\label{6-8}\\
\intertext{using standard microlocal analysis with logarithmic uncertainty principle and then in subsections 4.1--4.4 under assumption}
& \epsilon (h|\log h|)^{-\nu}\le \mu \le \epsilon h^{-\nu}\label{6-9}
\end{align}
applying microlocal analysis for $h$-pseudo-differential operators with respect to $x_2$ with operator-valued symbols -- operators in the auxiliary space $\bH = L^2(\bR_{x_1})$; I remind that in te case (\ref{6-9}) localization was done with respect to $\xi_2$ rather $x_1$.

Therefore in both cases (\ref{6-8}), (\ref{6-9}) in the redefined outer zone
\begin{equation}
\cZ_\out=\bigl\{{\bar\gamma}\le |x_1|\le {\bar\gamma}'_1=\epsilon {\bar\gamma}_1\bigr\}
\label{6-10}
\end{equation}
(with the small constant $\epsilon$) all these arguments remain true leading us eventually to the following statements:

\begin{proposition}\label{prop-6-1} Let conditions $(\ref{6-2})$ and $(\ref{6-4})_+$ be fulfilled. Let $\psi=\psi(x_2)$ be supported in
$B(0,{\frac 1 2})$ and let $\varphi=\varphi (\xi_2)$
be supported in the strip
\begin{equation}
\cY_\gamma=
\bigl\{\mu \gamma ^\nu \le |\xi_2|\le 2\mu \gamma^\nu\bigr\}
\label{6-11}
\end{equation}
with $C_1{\bar\gamma}\le \gamma \le \epsilon_1 {\bar\gamma}_1$
Then

\smallskip
\noindent
(i) As $\mu \le \epsilon h^{ -\nu}$ estimates
\begin{equation}
| F_{t\to h^{-1}\tau} \chi_T(t)\Gamma (Qe) |\le Ch^s
\label{6-12}
\end{equation}
and
\begin{equation}
\cR' = |\Gamma \bigl(Q e \bigr) - h^{-1}\int _{-\infty}^0 \Bigl( F_{t\to h^{-1}\tau} {\bar\chi}_T(t)\Gamma (Qe)\Bigr)\,d\tau |\le
C\mu^{-1} \gamma^{1-\nu}h^{-1}
\label{6-13}
\end{equation}
hold with $Qe=\varphi (hD_2) \bigl(e\psi\bigr)$, $e=e(x,y,\tau)$,
$|\tau|\le \epsilon$, $T\in [T_0,T_1]$, $T_0=Ch|\log h|$, $T_1=\epsilon \mu^{-1}\gamma^{-\nu}$;

\smallskip
\noindent
(ii) Moreover, under condition $(\ref{6-8})$ statement (i) holds with $Q=\psi_1\psi$, $\psi_1=\psi_1(x_1)$ supported in $\cZ_\gamma =
\{\gamma\le |x_1|\le 2\gamma\}$.
\end{proposition}

\begin{corollary} \label{cor-6-2} Let conditions $(\ref{6-2})$ and $(\ref{6-4})_+$ be fulfilled. Let $\psi=\psi(x_2)$
be supported in $B(0,{\frac 1 2})$ and $\varphi=\varphi(\xi_2)$
be supported in the outer zone defined in the terms of $\xi_2$
\begin{equation}
\cY_\out=
\bigl\{C_0\le |\xi_2|\le \epsilon \bigl(\mu h^\nu\bigr)^{-1/(\nu -1)}\bigr\}.
\label{6-14}
\end{equation}
Then

\smallskip
\noindent
(i) As $\mu \le \epsilon h^{ -\nu}$ estimate
\begin{equation}
\cR' \le C\mu^{-1/\nu} h^{-1}
\label{6-15}
\end{equation}
holds.

\smallskip
\noindent
(ii) Moreover, under condition $(\ref{6-8})$ statement (i) holds with $Q=\psi_1\psi$, $\psi_1=\psi_1(x_1)$ supported in $\cZ_\out$.
\end{corollary}

On the other hand, under condition  $(\ref{6-4})_-$ the whole zone $\cZ'=\cZ_\inn \cup \cZ_\out $ will be forbidden leading us to the following statement not having analogues in \cite{IRO6}:

\begin{proposition}\label{prop-6-3} Let conditions $(\ref{6-2})$ and $(\ref{6-4})_-$ be fulfilled. Let $\psi=\psi(x)$, $\psi_1=\psi_1(x_1)$ be supported in
$B(0,{\frac 1 2})$ and $\cZ'$ respectively and let
$\varphi =\varphi (\xi_2)$ be supported in the zone
\begin{equation}
\cY'=\{|\xi_2|\le \epsilon \bigl(\mu h^\nu\bigr)^{-1/(\nu -1)}\}.
\label{6-16}
\end{equation}
Then

\smallskip
\noindent
(i) $|Q e|\le Ch^s$ with $Qe=\varphi (hD_2) \bigl(e\psi\bigr)$, $e=e(x,y,\tau)$,
$|\tau|\le \epsilon$ as $\mu \le \epsilon h^{-\nu}$;

\smallskip
\noindent
(ii) Moreover, under condition $(\ref{6-8})$ statement (i) holds with $Q=\psi_1\psi$, $\psi_1=\psi_1(x_1)$ supported in $\cZ'$.
\end{proposition}

Therefore as $\mu \le \epsilon h^{ -\nu}$ and condition $(\ref{6-4})_+$ is fulfilled one needs to discuss the contribution of the inner zone
$\cZ_\inn =\{|x_1|\le {\bar\gamma} \}$ or equivalently
$\cY_\inn = \{|\xi_2|\le C_0\}$\,\footnote{\label{foot-3} These two definitions are essentially equivalent under condition (\ref{6-8}) but in the case (\ref{6-9}) one needs always use definition in the frames of $\xi_2$.} to the remainder estimate.
Furthermore one needs to consider the contribution of the previously forbidden zone
$\cZ''=\{|x_1|\ge {\bar\gamma}'_1\}$ or equivalently
$\cY''=\{|\xi_2|\ge \epsilon\bigl(\mu h^\nu\bigr)^{-1/(\nu -1)}\}$\,$^{\ref{foot-3}}$
to the remainder estimate.

The inner zone is analyzed exactly as in section 2 and subsections 4.1--4.4 of \cite{IRO6} leading us eventually to

\begin{proposition}\label{prop-6-4-pm} Let conditions $(\ref{6-2})$ and $(\ref{6-4})_+$ be fulfilled. Let $\psi=\psi(x_2)$ and $\psi_1=\psi_1(x_1)$ be supported in
$B(0,{\frac 1 2})$ and $\cZ_\inn$ respectively and let
$\varphi =\varphi (\xi_2)$ be supported in $\cY_\inn= \{ |\xi_2|\le C_0\}$. Then
all the results of section 2 and subsections 4.1--4.4 of {\rm\cite{IRO6}} remain true; in particular

\smallskip
\noindent
(i) As $\mu \le Ch^{\delta -\nu}$
\begin{equation}
\cR'' \Def |\Gamma \bigl(Q e \bigr) - h^{-1}\sum _j \int _{-\infty}^0 \Bigl( F_{t\to h^{-1}\tau} {\bar\chi}_{T_j}(t)\Gamma (Q_je)\Bigr)\,d\tau |\le C\mu^{-1/\nu} h^{-1}
\label{6-17}
\end{equation}
with $Qe=\varphi (hD_2) \bigl(e\psi\bigr)$, $e=e(x,y,\tau)$, $Q=\sum_j Q_j$ and
$|\tau|\le \epsilon$ where partition $Q_j$ and $Ch|\log h|\le T_j$ are defined following formula $(3.28)$ in {\rm \cite{IRO6}};

\smallskip
\noindent
(ii) Moreover, under nondegeneracy condition
\begin{pdeq}\label{6-18}\end{pdeq}
\begin{equation}
\sum_{1\le k\le m}|\partial_{x_2}^k \bigl({\frac W f}\bigr)|\ge \epsilon_0.
\label{6-18-m}
\tag*{$(\theequation)_m$}
\end{equation}
$\cR''$ does not exceed $C\mu^{-1/\nu} h^{-1}$ as $\mu \le \epsilon h^{-\nu}$;

\smallskip
\noindent
(iii) On the other hand, in the general case $\cR''$ does not exceed
$C\mu^{-1/\nu} h^{-1}+Ch^{-\delta}$ as $\mu \le \epsilon h^{-\nu}$;

\smallskip
\noindent
(iv) Furthermore, under condition $(\ref{6-8})$ all statements (i)--(iii) hold with $Q=\psi_1\psi$.
\end{proposition}

\begin{remark}\label{rem-6-5} In frames of proposition \ref{prop-6-4-pm} estimate
(\ref{6-12}) holds for $Q=Q_m$ and $T\in [T_m, T'_m]$ with $T'_m$ defined by (\ref{IRO6-2-98}) (it was denoted by $T_1$ then).
\end{remark}

\subsection{Estimates. II}

To investigate zone $\cZ''$ I will apply the theory of operators with operator-valued symbols. However, as $\mu \le \epsilon (h|\log h|)^{-\nu}$ one can apply a usual microlocal analysis with logarithmic uncertainty principle.

So, let us consider $A$ as $h$-pseudo-differential operator $\cA (x_2, h D_2)$ with operator-valued symbol $\cA(x_2,\xi_2)$. However, before doing this one can assume without any loss of the generality that $g^{11}=1$, $g^{12}=0$ and therefore
\begin{multline}
\cA(x_2,\xi_2)=
{\frac 1 2}\biggl(h^2D_1^2 + \sigma^2 (x) \bigl(\xi_2 -\mu V_2(x)\bigr)^2 - (2l+1)\mu h F -W(x)\biggr), \\
V_2=\phi (x){\frac 1 \nu}x_1^\nu
\label{6-19}
\end{multline}
with $\phi (x)=1$ as $x_1=0$; then $f=\sigma\phi$.

Further, for given $x_2$ by change of variable $x_1$ one can transform $\cA$ unitarily to the similar operator with $\phi =1$ and with
\begin{equation}
\sigma =1 \qquad \text{as\;} x_1=0;
\label{6-20}
\end{equation}
but this new operator is multiplied from the left and the right by $\alpha (x)$.
So operator $\cA(x_2,\xi_2)$ is unitary equivalent to
\begin{multline}
\cA' (x_2,\xi_2)= \\
{\frac 1 2}\alpha (x)
\biggl(h^2D_1^2 + \sigma^2 (x) \Bigl(\xi_2 -\mu {\frac 1 \nu}x_1^\nu \Bigr)^2 -(2l+1)\mu h \sigma (x)x_1^{\nu -1}-W_0(x)\biggr)
\alpha (x).
\label{6-21}
\end{multline}
Note that $W_0=W/f$ as $x_1=0$ and thus \emph{conditions \ \ref{6-4-pm}, $(\ref{6-5})$ and \ref{6-18-m} are reformulated in terms of $W_0$\/} obviously.

Proposition \ref{prop-A-3}(ii) of Appendix A implies that under condition \ref{6-4-pm} zone
$\cY''\setminus \cY''_0=\{\epsilon (\mu h^\nu)^{-1/(\nu -1)} \le |\xi_2|\le
2C(\mu h^\nu)^{-1/(\nu -1)}\}$ is microhyperbolic with respect to $\xi_2$ and thus one can extend $\cY'$ to zone
${\bar\cY}' \Def \{ |\xi_2|\le 2C(\mu h^\nu)^{-1/(\nu -1)}\}$ resulting in the following statement:

\begin{proposition}\label{prop-6-6} Let conditions $(\ref{6-2})$ and\ref{6-4-pm} be fulfilled. Then estimate $\cR'\le C$ holds as $\cR'$ is defined by $(\ref{6-13})$ with $Qe=\varphi (hD_2)(e\psi)$, $\varphi $ supported in the zone $\cY''\setminus \cY''_0$, $T\in [T_0,T_1]$, $T_0=Ch|\log h|$,
$T_1=\epsilon (\mu h^\nu )^{-1/(\nu-1)}$, $\mu\le \epsilon h^{-\nu}$.
\end{proposition}

Furthermore,  proposition \ref{prop-A-3}(i) implies that under condition \ref{6-4-pm} zone
$\cY''_0=\{|\xi_2|\ge C(\mu h)^{-1/(\nu-1)}\}$ is forbidden on energy levels $|\tau|\le \epsilon$ as long as $\mu \le \epsilon h^{-\nu}$ is forbidden; namely
\begin{equation}
| F_{t\to h^{-1}\tau }{\bar\chi}_T(t) (Qu)(x,y,t) |\le CT h^s
\qquad \forall \tau:|\tau|\le\epsilon
\label{6-22}
\end{equation}
as $Q\psi =\varphi (hD_2)(u\psi)$ with $\varphi $ supported in the zone $\cY''_0$ and therefore its contribution to the remainder $\cR'$ defined by (\ref{6-13}) is negligible as well:

\begin{proposition}\label{prop-6-7} Let conditions $(\ref{6-2})$ and \ref{6-4-pm} be fulfilled. Then estimate
$\cR'\le Ch^s$ holds as $\cR'$ is defined by $(\ref{6-13})$ with
$Qe=\varphi (hD_2)(e\psi)$, $\varphi $ supported in the zone $\cY''_0$,
$T\ge T_0=Ch|\log h|$, $\mu\le \epsilon h^{-\nu}$.
\end{proposition}

The analysis of all zones under condition (\ref{6-5}) will be done in subsection 6.7.

\subsection{Calculations. I}

In this subsection I will change partition: instead of $\cZ'$ and $\cZ''$ I will consider ${\bar\cY}'$ and $\cY''_0$ obtained if I redefine
${\bar\gamma}_1= C(\mu h)^{-1/(\nu -1)}$; respectively change definitions and notations of zones $\cY_\out$, $\cZ_\out$, $\cZ'$, $\cZ''$.

After estimates were derived in two previous subsections under assumption
$Ch^{-1}\le \mu \le \epsilon h^{-\nu}$ and condition \ref{6-4-pm} calculations in zone ${\bar\cY}'$ are done exactly as in section 3 and subsection 4.4 of \cite{IRO6}.

On the other hand, calculations in zone $\cY''_0$ as $\mu \le \epsilon h^{-\nu}$ are rather obvious under assumptions
$Ch^{-1}\le \mu \le \epsilon h^{-\nu}$ and \ref{6-4-pm}. Therefore I arrive to the intermediate estimate
\begin{equation}
|\int \Bigl( \bigl(\varphi (hD_2)e\bigr)(x,x,0)\ - (2\pi h)^{-1}
\int {\rm e}(x_1,x_1;x_2,\xi_2,0)\varphi(\xi_2)\,d\xi_2\Bigr)
\psi _2(x_2)\, dx|\le R
\label{6-23}
\end{equation}
where $R$ is an estimate already derived in the corresponding conditions (also see below) and $\varphi\in C_0^\infty (-\epsilon',\epsilon')$ with sufficiently small constant $\epsilon'$.

Then the same estimate holds with $\psi(x_2)$ replaced by $\psi (x)$ such that
$\psi (x)=\psi _2(x_2)$ as $|x_1|\le C_1\epsilon'$ because this transition leads to a negligible error. I take $\psi$ also satisfying $\psi(x)=0$ as $|x_1|\ge 2C_1\epsilon'$. Then in the latter estimate I can replace $\varphi $ by 1. Really, then the error would be
\begin{equation}
|\int \Bigl( \bigl((1-\varphi (hD_2))e\bigr)(x,x,0)\ - (2\pi h)^{-1}
\int {\rm e}(x_1,x_1;x_2,\xi_2,0)(1-\varphi(\xi_2))\,d\xi_2\Bigr)
\psi (x)\, dx|
\label{6-24}
\end{equation}
and replacing $\psi$ by $\psi'$ equal to $\psi$ as $|x_1|\ge 2C_2^{-1}\epsilon'$
and equal to 0 as $|x_1|\le C_2^{-1}\epsilon'$ leads to a negligible error. However, to expression (\ref{6-24}) modified this way one can apply the theory of operators with non-degenerating magnetic field and then to estimate expression (\ref{6-24}) by $C$.

Thus I derived (\ref{6-23}) with $\varphi$ replaced by 1 and $\psi_2(x_2)$ replaced by some ``special'' function $\psi(x)$. Then due to rescaling arguments like in subsubsection 6.2.2 the same estimate holds for a general function $\psi(x)$ supported in $\{|x_1|\le 2C_1\epsilon'\}$. Thus I arrive to

\begin{proposition}\label{prop-6-8} Let conditions $(\ref{6-2})$ and $(\ref{6-4})_+$ be fulfilled. Then 

\smallskip
\noindent
(i) As \underline{either} $\mu \le h^{\delta-\nu}$ \underline{or} condition \ref{6-18-m} is fulfilled and $\mu \le \epsilon h^{-\nu}$ the following estimate holds
\begin{equation}
\cR_I\Def |\int \Bigl( e(x,x,0)\ - (2\pi h)^{-1}
\int {\rm e}(x_1,x_1;x_2,\xi_2,0)\,d\xi_2\Bigr)
\psi (x)\, dx|\le C\mu^{-{\frac 1 \nu }}h^{-1}
\label{6-25}
\end{equation}
where here and below ${\rm e}(x_1,y_1; x_2,\xi_2,\tau)$ is the Schwartz kernel of the spectral projector of operator $\cA(x_2,\xi_2)$ defined by $(\ref{6-19})$ and $\delta >0$ is an arbitrarily small exponent;

\smallskip
\noindent
(ii) In the general case with $\mu \le \epsilon h^{-\nu}$ estimate
\begin{equation}
\cR_I \le C\mu^{-{\frac 1 \nu }}h^{-1} +Ch^{-\delta}
\label{6-26}
\end{equation}
holds.
\end{proposition}

I remind that in both statements of proposition \ref{prop-6-8} the principal part of asymptotics has magnitude $\asymp \mu h^{-1}$ (as $\mu \ge h^{-1}$).

On the other hand, under condition
$(\ref{6-4})_-$ zone $\cY'$ becomes forbidden and thus I arrive to

\begin{proposition}\label{prop-6-9} Let conditions $(\ref{6-2})$ and $(\ref{6-4})_-$ be fulfilled and $l\ge 1$. Then for
$Ch^{-1}\le \mu \le \epsilon h^{-\nu}$ estimate $\cR_I\le C$ holds while the principal part of asymptotics has magnitude $\asymp \mu h^{-1}$.
\end{proposition}

\subsection{Calculations. II}

Transition to the auxiliary operator $\cA_0$ without increasing error estimates could be done easily in zone $\cY_\out$ exactly as it was done in the proof of propositions \ref{IRO6-prop-3-3} and \ref{IRO6-prop-3-4} while arguments of \ref{IRO6-prop-3-8} etc work in zone $\cY_\inn$. 

On the other hand, this transition in zone $\cY''_0$ is obvious under condition \ref{6-4-pm}, and I arrive to two theorems below as $\mu \le h^{-\nu}|\log h|^{-K}$ and function $\psi$ is ``special'' in the sense of the previous subsection. Then the same arguments as there extend theorem to general $\psi$.

Furthermore, under condition \ref{6-4-pm} the case
$h^{-\nu}|\log h|^{-K}\le \mu \le \epsilon h^{-\nu}$ is analyzed exactly as in section 4 of Part I leading to the extension of these theorems to
$\mu \le \epsilon h^{-\nu}$:

\begin{theorem}\label{thm-6-10} Let conditions $(\ref{6-2})$ and $(\ref{6-4})_+$ be fulfilled. Then

\smallskip
\noindent
(i) As \underline{either} $\mu \le h^{\delta-\nu}$ \underline{or} condition \ref{6-18-m} is fulfilled and $\mu \le \epsilon h^{-\nu}$ 
\begin{align}
&\cR^*\Def |\int \Bigl( e(x,x,0)- {\tilde\cE}^\MW (x,0)\Bigr)\psi (x)\,dx - \int \cE^\MW_\corr (x_2,0) \psi(0,x_2)\,dx_2|
\label{6-27}\\
\intertext{does not exceed $C\mu^{-1/\nu}h^{-1}$ where}
&\cE^\MW_\corr (x,\tau)\Def (2\pi h)^{-1}
\int {\rm e}_0(x_1, x_1;x_2, \xi_2, \tau, \hbar) \,dx_1 d\xi_2-
\int {\tilde \cE}^\MW_0 (x,\tau)\,dx_1,
\label{6-28}
\end{align}
$\cE^\MW$ is Magnetic Weyl approximation\footnote{\label{foot-4} See e.g. $(\ref{IRO6-0-8})$.} and here and below ${\rm e}_0(x_1,y_1; x_2,\xi_2,\tau)$ is the Schwartz kernel of the spectral projector of operator $\cA_0(x_2,\xi_2)$ defined by $(\ref{6-19})$ and with $\alpha, \phi, \sigma, W$ restricted to $\{x_1=0\}$ and $\cE^\MW_0$ is Magnetic Weyl approximation for this operator.

\smallskip
\noindent
(ii) In the general case with $\mu \le \epsilon h^{-\nu}$ estimate
$\cR^*\le C\mu^{-1/\nu}+Ch^{-\delta}$ holds.
\end{theorem}

\begin{theorem}\label{thm-6-11} Let conditions $(\ref{6-2})$ and $(\ref{6-4})_-$ be fulfilled and $l\ge 1$. Then as
$Ch^{-1}\le \mu \le \epsilon h^{-\nu}$ estimate $\cR^*\le C$ holds while the principal part of asymptotics has magnitude $\asymp \mu h^{-1}$.
\end{theorem}

\begin{remark}\label{rem-6-12} 
Obviously the same approximate expressions  (\ref{IRO6-3-52}), \ref{IRO6-3-52-*}, \ref{IRO6-3-52-**} hold for the part of $\cE^\MW_\corr$ ``associated'' with $\cY_\inn$;
\end{remark}

\subsection{Estimates under condition $(\ref{6-5})$}

I start from the remainder estimate in zone ${\bar\cY}'$ which is trivial:

\begin{proposition}\label{prop-6-13} Let conditions $(\ref{6-2})$, $(\ref{6-20})$ and $(\ref{6-5})$ be fulfilled. Then

\smallskip
\noindent
(i) Estimate $(\ref{6-13})$ holds with $Qe=\varphi (hD_2)(e\psi)$, $\varphi $ supported in the strip $\cY_\gamma$ with the same restrictions to $\gamma$ and
the same $T_0$, $T_1$ as in proposition \ref{6-1}(i);

\smallskip
\noindent
(ii) Furthermore, the same estimate holds as $\varphi $ is supported in zone $\cY_\inn$ and $\gamma={\bar\gamma}_0=\mu^{-1/\nu}$;

\smallskip
\noindent
(iii) Therefore $\cR'$ defined by $(\ref{6-13})$ does not exceed  $C\mu^{-1/\nu}h^{-1}$ as $\varphi$ is supported in zone ${\bar\cY}'$ and $T=T_0$.
\end{proposition}

Now let us analyze zone $\cY''_0$ under condition $(\ref{6-5})$:

\begin{proposition}\label{prop-6-14} Let conditions $(\ref{6-2})$, $(\ref{6-20})$ and $(\ref{6-5})$ be fulfilled. Then estimate
$\cR'\le C$ holds as $\cR'$ is defined by $(\ref{6-13})$ with
$Qe=\varphi (hD_2)(e\psi)$, $\varphi $ supported in the zone $\cY''_0$.
\end{proposition}

\begin{proof} (i) Let us note first that estimate
\begin{equation}
| F_{t\to h^{-1}\tau }\bigl({\bar\chi}_{T_1}(t)- {\bar\chi}_{\bar T}(t)\bigr) (Qu)(x,y,t) |\le Ch^s
\qquad \forall \tau:|\tau|\le\epsilon
\label{6-29}
\end{equation}
holds with $T_1= \epsilon \mu^{-1} \gamma^{-\nu}$, ${\bar T}=Ch|\log h|$
as $Qu=\varphi (hD_2)(u\psi)$, $\varphi $ supported in the strip
$\cY_{(\gamma)}=\{\mu \gamma^\nu \le |\xi_2|\le 2 \gamma^\nu\}$ with
$\gamma \ge C{\bar\gamma}_1$.

Really, us consider a partial trace $\Gamma' (Qu)$ (with respect to $x_1$). Due to proposition \ref{prop-A-3} the propagation speed with with respect to $x_2$ does not exceed $C|\xi_2|^{-1}\asymp C(\mu \gamma^\nu)^{-1}$ and the propagation speed with respect to $\xi_2$ does not exceed $C$\,\footnote{\label{foot-5} Under some assumptions this would be equivalent to the estimate of the the average propagation speed with respect to $x_1$ of $Qu$ by
$C\gamma (\mu \gamma^\nu)^{-1}$; further one can estimate
average propagation speed with respect to $x_2$ of $Qu$ by
$C(\mu \gamma^\nu)^{-1}$
as well.}; moreover, under condition (\ref{6-5}) this propagation speed with respect to $\xi_2$  is greater than $\epsilon$.

On the other hand, an obvious estimate
\begin{equation}
| F_{t\to h^{-1}\tau } {\bar\chi}_{T_0}(t) \Gamma(Qu)(t) |\le
C\mu\gamma^\nu h^{-1}\times T_0= C\mu \gamma^\nu |\log h|
\label{6-30}
\end{equation}
holds where the first factor is $\mu_\eff h_\eff^{-1}\gamma^{-1}$; furthermore, due to (\ref{6-29}) this estimate holds for the left-hand expression with $T_0$ replaced by $T_1$.

Therefore the contribution of the strip $\cY_\gamma$ to the remainder estimate does not exceed
\begin{equation}
C\mu \gamma^\nu |\log h| \times T_1^{-1}= C|\log h|
\label{6-31}
\end{equation}
and therefore the total contribution of $\cY''_0$ to the remainder estimate does not exceed
$C|\log h| \int \gamma^{-1}\,d\gamma \asymp C|\log h|^2$.

This estimate is as good as I need for $\mu \le Ch^{-\nu}|\log h|^{-2\nu}$. However for $Ch^{-\nu}|\log h|^{-2\nu}\le \mu \le \epsilon h^{-\nu}$ I would like to improve it getting rid of two logarithmic factors.

\medskip
\noindent
(ii) Getting rid off one of them is easy: rescaling $t \mapsto t/T$, $(x_j-y_j)\mapsto (x_j-y_j)/T$, $\mu \mapsto \mu T$, $h\mapsto h/T$
estimates for Schr\"odinger operator with strong non-degenerate magnetic field
\cite{Ivr1}, section 6 (with arbitrary parameters $\mu$ and $h$ such that 
$\mu h \ge C$) I arrive to two following inequalities
\begin{align}
& |F_{t\to h^{-1}\tau}\chi_T(t) \Gamma (Qu) |\le C\mu \bigl({\frac h T} \bigr)^s
\label{6-32}\\
& |F_{t\to h^{-1}\tau}{\bar \chi}_T(t) \Gamma (Qu) |\le C\mu \label{6-33}
\end{align}
as $h\le T\le 1$, $|\tau|\le \epsilon$ under condition
$|W|+|\nabla W|\ge \epsilon_0$. Then using our standard scaling
$x_1\mapsto x_1/\gamma$, $x_2\mapsto (x_2-y_2)/\gamma$,
$\mu \mapsto \mu_\eff =\mu \gamma^\nu$, $h\mapsto h_\eff=h/\gamma$ and 
$T\mapsto T/\gamma$ I arrive to estimate (\ref{6-30}) without logarithmic factor
\begin{equation}
| F_{t\to h^{-1}\tau } {\bar\chi}_T(t) \Gamma(Qu)(t) |\le
C\mu \gamma^\nu
\tag*{$(\ref*{6-30})^*$}\label{6-30-*}
\end{equation}
as $|\tau|\le\epsilon$, $T/\gamma \le \epsilon \mu \gamma^\nu \iff T\le T'_1=\epsilon \mu \gamma^{\nu +1} $. Further,  and  to
$(\ref{6-29})$ this estimate holds as $h\le T\le T_1=\epsilon \mu \gamma^\nu $ provided $T'_1\ge Ch$ i.e. $\gamma \ge {\bar \gamma}_1$.

Then the contribution of the strip $\cY_\gamma$ to the remainder $\cR'$ is $C$ and therefore the total estimate is $C|\log h|$.

\medskip
\noindent
(iii) To get rid off the second logarithmic factor I need to further increase $T_1$ in the previous arguments and for this purpose I need for each $\gamma$ to make $x_2$-partition of $\cY_\gamma$ of the size
\begin{equation}
\ell = \epsilon |V(0,x_2)| +{\bar \ell},\qquad {\bar\ell}\ge C\gamma.
\label{6-34}
\end{equation}

Consider first elements $\cU_{\gamma,\ell}$ with $\ell\ge C{\bar\ell}$. For every such element on levels $\tau$ with $|\tau|\le \epsilon \ell$ after rescaling
\begin{equation}
x_2\mapsto x_2 \ell^{-1},\quad h\mapsto h' =h \ell ^{-\frac 3 2},\quad
t\mapsto t \ell^{-1},\quad
\mu\mapsto \mu'=\mu \ell^{\frac 1 2}
\label{6-35}
\end{equation}
I am in the elliptic situation.

Therefore contribution of each such element to the remainder estimate does not exceed $C\mu '(h')^s$ and therefore the total  contribution of such elements is negligible as ${\bar\ell}= h^\delta$.

So I need to consider only elements $\cU'_\gamma=\cU_{\gamma,\ell}$ with
$\ell \asymp {\bar\ell}=h^\delta$. For such elements after rescaling (\ref{6-35}) I can apply estimate \ref{6-30-*}; then scaling back I get the same estimate \ref{6-30-*} again but with $Q=\psi' (x_2)\varphi (hD_2)$
supported in $\cU'_\gamma$, $|\tau|\le \epsilon \ell$ and
$Ch|\log h|\ell^{-1}\le T\le T_1=\epsilon \mu \gamma^{\nu+1}$\,\footnote{\label{foot-6} It is consistent with the fact that support of $\psi'$ is of the length $\ell$ but now ${\bar T}=Ch|\log h|/\ell$. }. Furthermore, applying (\ref{6-29}) I can increase $T_1$ to 
$\epsilon \mu \gamma^\nu$.

So far I gained nothing: the estimate I  proved alone would bring me the same final remainder estimate $C|\log h|$ as before but now I can further increase $T_1$  and thus reduce the remainder estimate.

Namely, let us consider propagation in the time direction in which $|\xi_2|$ increases. If only propagation with respect to $\xi_2$ was considered, until time $\epsilon_3 \mu$ it would be confined to zone
\begin{equation*}
\bigl\{\epsilon_0 \le |\xi_2|\bigl(\mu \gamma^\nu + |t|\bigr)^{-1}\le C\bigr\} \subset
\bigl\{{\frac 1 2}\mu \gamma^\nu \le |\xi_2|\le \epsilon _1\mu\bigr\}
\end{equation*}
and thus to
$\{|x_1|\le \epsilon_3\}$. 

However let us note that the propagation speed with respect to $x_2$ does not exceed $C\ell/|\xi_2|$ as $\ell\ge C|V|+{\bar\ell}$. Therefore one can prove easily that propagation, which started in the zone
$\{|x_2|\le {\frac 1 2}, |V|\le h^\delta\}$ as I have assumed,
until time $T^*_1=\mu \gamma^\nu h^{-\delta_1}$ is confined to a bit larger zone
$\bigl\{|x_2|\le {\frac 3 4}, |V|\le h^{\delta/2}\bigr\}$ of the same type.

Therefore estimate \ref{6-30-*} holds with $Ch ^{1-2\delta} \le T\le T^*_1$. Then due to the Tauberian approach contribution of each partition element $\cU'_\gamma$ to the remainder estimate does not exceed $
C\mu \gamma^\nu T^{* -1}_1= Ch^{\delta_1}$ and the contribution of the whole strip $\cY_\gamma$ does not exceed $Ch^{\delta_1}$ as well and of the whole zone $\cY''_0$ does not exceed $Ch^{\delta_2}$.

Clearly, at some moment I increased slightly $T_0$ but after summation over partition was done I can (using negligibility of the trace on
$[Ch|\log h|, h^{1-\delta}]$ time interval on energy levels
$|\tau|\le \epsilon$) return to original ${\bar T}$.
\end{proof}

\subsection{Calculations under condition $(\ref{6-5})$}

Calculations in zone ${\bar\cY}'$ are exactly as in \cite{IRO6}. However one should be more careful with calculations in zone $\cY''_0$.

Let me remind that according to subsection 6.2 \cite{Ivr1} in the nondegenerate case with $\mu h\ge C$ the operator in question is reduced to one-dimensional $\mu^{-1}h$-pdo  $B(x_2,\mu^{-1}hD_2, h^2)$\,\footnote {\label{foot-7} Where $x_2$ is not our original $x_2$.} with the ``main symbol'' $B(x_2,\xi_2, 0)= W \circ \Psi$ and therefore the contribution of the partition element to the final answer will be given as in subsection 6.6 by magnetic Weyl expression
$\int \cE^\MW (x,0)\psi (x)\, dx$ plus correction terms
$ \mu h ^{1+2m} \int \varkappa_{l,m} (x)\psi (x) \, dx$, $m=0,1,\dots$.

After rescaling $\mu\mapsto \mu \gamma^\nu$, $h\mapsto h/\gamma$, $dx\mapsto \gamma^{-2}dx$ these terms are transformed into
\begin{equation}
\mu h^{1+2m} \int \varkappa_{l,m} (x,\gamma )\psi (x) \gamma^{\nu -2m -3} \,dx
\label{6-36}
\end{equation}
integrated over zone $\{{\bar \gamma}_1\le \gamma \le \epsilon\}$.

One can see easily that if there was an extra factor $\gamma$ one would be able to rewrite this expression (\ref{6-36}) modulo $O(1)$ into the similar expression with integration over $\{\gamma \le \epsilon\}$ as
$2m+2 < \nu$\,\footnote{\label{foot-8} thus resulting in exactly expression $\kappa_{l,m} \mu h^{1+2m}$ as in non-degenerate case.}
or to simply skip it as $2m+2 > \nu$ or to get a term which is
$O\bigl(\mu h^\nu |\log h|\bigr)$ as $2m+2=\nu$. To gain this extra factor one needs to consider the difference of expressions $\int e(x,x,0)\psi (x)\,dx$ for two operators with $g^{jk}(x)$, $f(x)$, $V(x)$ coinciding as $x_1=0$. As this second operator it is natural to pick up the simplest one i.e.
\begin{equation}
A_0= {\frac 1 2}\biggl(h^2D_1^2 + \bigl(hD_2-\mu x_1^\nu/\nu\bigr)^2 - (2l+1)\mu h x_1^{\nu-1} - W(x_2)\biggr).
\label{6-37}
\end{equation}
Therefore I arrive to

\begin{proposition}\label{prop-6-15} Under condition $(\ref{6-5})$ estimate
\begin{multline}
|\int \biggl( e(x,x, 0) - e_0(x,x,0) -\cE^\MW (x,0) + \cE^\MW_0 (x,0)\biggr)
\psi (x)\,dx - \\
\sum \kappa'_{l,m} \mu h^{1+2m}| \le C\mu^{-{1/\nu}}h^{-1}
\label{6-38}
\end{multline}
holds as $\mu \le h^{-\nu}|\log h|^{-K}$ where $e_0$ and $\cE^\MW_0$ are defined
for operator $A_0$.
\end{proposition}

\begin{claim}\label{6-39}
Now in what follows I can consider operator $A_0$ instead of $A$.
\end{claim}

Then I can apply the standard method of successive approximations with unperturbed operator $\cA (y_2,hD_2)$ and plug the results of successive approximations into expression
\begin{equation}
h^{-1}\int_{-\epsilon}^0 \biggl( F_{t\to h^{-1}\tau} {\bar\chi}_T(t) \Gamma (Qu)\biggr)\,d\tau
\label{6-40}
\end{equation}
which calculates exactly contribution of the ``problematic'' eigenvalue $\lambda_l$ of the corresponding one-dimensional Schr\"odinger operator; I remind that $T={\bar T}=Ch|\log h|$.

Thus while the main part of asymptotics is estimated by
$C\mu h^{-2}\gamma ^\nu T = C\mu h^{-1}\gamma^\nu |\log h|$, each next term seemingly acquires factor
\begin{equation}
Ch^{-1}\bigl(\mu h \gamma^{\nu-1}\bigr)^{1/2}T^2 \asymp
Ch \bigl(\mu h \gamma^{\nu-1}\bigr)^{1/2}|\log h|^2;
\label{6-41}
\end{equation}
since the propagation speed with respect to $x_2$ is estimated by
$C_0(\mu h\gamma^{\nu -1})^{1/2}$ such factor could be larger than 1.

In fact however, $C_0(\mu h\gamma^{\nu -1})^{1/2}$ is the estimate for the instant propagation speed only. Using instead  the mentioned reduction to a one-dimensional $\mu^{-1}h$-pdo one can find that the propagation speed with respect to $x_2$ is estimated by $C_0\mu ^{-1}$ if magnetic field is non-degenerate and then in the canonical coordinates for time $T={\bar T}$ the shift of $(x'_2,\xi'_2)$ will be estimated by $C_0\bigl(\mu^{-}h |\log h|\bigr)^{1/2}$ which is the smallest distance allowed by the logarithmic uncertainty principle\footnote{\label{foot-9} Since $\mu^{-1}h$-Fourier Integral Operators are involved later one needs the same distance in each $(x,\xi)$ direction.} and this would persist if one returns back to the original 
$(x_2, \mu^{-1}\xi_2)$; so one would be able to estimate $(x_2-y_2)$ on the time interval in question by
$C_0\bigl(\mu^{-1}h |\log h|\bigr)^{1/2}$.

In the degenerate case described here one must replace $\mu$, $h$ by
$\mu \gamma^\nu$, $h/\gamma$ respectively and then multiply by $\gamma$ thus producing final estimate for $|x_2-y_2|$
\begin{equation}
\varrho\Def C\bigl(\mu ^{-1}h \gamma ^{1-\nu} |\log h|\bigr)^{1/2} \asymp
Ch {\bar \gamma}_1^{{\frac 1 2}(\nu -1)}\gamma ^{-{\frac 1 2}(\nu -1)}
|\log h|^{\frac 1 2}
\label{6-42}
\end{equation}
and therefore each next term acquires factor $\varrho |\log h|$. Then $m$-th term of the final answer is estimated by
\begin{equation}
C\mu h^{-1}\varrho^{m-1} |\log h|^K \asymp C\mu h^{m-2}
\gamma^{\nu - {\frac 1 2}(\nu -1)(m-1)}
{\bar\gamma}_1^{{\frac 1 2}(\nu-1)(m-1)}|\log h|^K.
\label{6-43}
\end{equation}
After integration over $\gamma^{-1}\,d\gamma$ with
${\bar\gamma}_1\le \gamma \le \epsilon$ expression (\ref{6-43}) results in
$C\mu h^{m-2}{\bar\gamma}_1^\nu|\log h|^K$ as
$\nu - {\frac 1 2}(\nu -1)(m-1)\le 0$ or in
$C(\mu ^{-1}h)^{(m-3)/2}|\log h|^K$ otherwise.
One can check easily that in either case the answer is $O(|\log h|^K)$ as
$m\ge 3$ and only terms with $m=1,2$ should be considered more carefully under condition (\ref{6-8}).

\smallskip
On the other hand, the main term appears as (\ref{6-40}) with $u$ replaced by ${\bar u}$ and modulo negligible one can rewrite it with any $T\ge {\bar T}$, in particular with $T=\infty$ which leads to
\begin{equation}
(2\pi h)^{-1} \int
{\rm e}_0(x_1,x_1,0; x_2,\xi_2)\psi (x_1)\varphi (\xi_2) \,dx_2\,d\xi_2
\label{6-44}
\end{equation}
where I remind that ${\rm e}_0(x_1,y_1,0; x_2,\xi_2)$ is the Schwartz kernel of the spectral projector of one-dimensional Schr\"odinger operator $\cA_0(x_2,\xi_2)$.

\smallskip
Let us consider terms with $m=2$ i.e. expression (\ref{6-40}) with $u$ replaced by ${\bar u}_1$; similarly to analysis of (i) one can estimate contribution of $O\bigl((x_2-y_2)^2\bigr)$ terms in the perturbation
$\cA (x_2,hD_2) - \cA (y_2,hD_2)$ by $C|\log h|^K$. Therefore one should consider only
$\cA (x_2,hD_2) - \cA (y_2,hD_2)= (x_2-y_2)B_1 (y_2)$ in which case
${\bar u}_1$ is defined by (3.23) without the last term since $B_1$ commutes with $(x_2-y_2)$:
\begin{equation}
u \mapsto {\bar u}_1=-ih\sum_{\varsigma=\pm }\varsigma {\bar G}^\varsigma B_1
{\bar G}^\varsigma [{\bar A}, x_2-y_2] {\bar G}^\varsigma \delta(t)\delta (x_2-y_2)\delta(x_1-y_1).
\label{6-45}
\end{equation}

One needs to multiply this by $h^{-1}\psi$, integrate with respect to $\tau$ and apply $\Gamma$ to it. Obviously since
for odd $\nu$ operators ${\bar G}^\varsigma$ and $[{\bar A}, x_2-y_2]$ are even and odd respectively as $x_1\mapsto -x_1$, $\xi_2\mapsto -\xi_2$ the answer would be 0 if $\psi$ is even with respect to $x_1$.

To cover the case of even $\nu$ and general $\psi$ let us note that $B_1$ commutes with ${\bar G}^\varsigma$ considered as operators in the auxiliary space $L^2(\bR^1_{x_1})$. Then if ${\bar G}^\varsigma$ commuted with $\psi$, taking trace and integrating with respect to $\tau$ would result in
\begin{equation*}
\const \cdot \partial_{\xi_2} B_1 \sum_{\varsigma=\pm }
\varsigma \Tr \bigl( {\bar G}^\varsigma \psi\bigr)
\end{equation*}
which after integration over $\xi_2$ results in $0$.

However ${\bar G}^\varsigma$ does not commute with $\psi$, so instead of 0 one gets
\begin{equation*}
\const \cdot B_1 \sum_{\varsigma=\pm }\Tr
\varsigma\Bigl({\bar G}^\varsigma
\bigl(\partial_{\xi_2} {\bar G}^\varsigma\bigr) \Bigl({\bar G}^\varsigma
[{\bar A},\psi]\Bigr)
\end{equation*}
and to this expression one can apply the same type of transformations and calculations as in the proof of proposition \ref{prop-6-15} resulting in
the expressin $\sum_m \kappa_{l,m} \mu h^{1+2m}$ where coefficients $\kappa_{l,m}$ are changed as needed.

Therefore combining with the results for zone ${\bar\cY}'$ I arrive to

\begin{proposition}\label{prop-6-16} For a model operator
\begin{multline}
|\int \biggl(e_0(x,x,0) - (2\pi h)^{-1}\int
{\rm e}_0(x_1,x_1,0; x_2,\xi_2)\,d\xi_2\biggr)\psi (x)\,dx -\\
\sum \kappa_m \mu h^{1+2m}|\le C\mu^{-1/\nu}h^{-1}
\label{6-46}
\end{multline}
as $\mu \le Ch^{-\nu}|\log h|^{-K}$.
\end{proposition}

Further, combining this with proposition \ref{prop-6-14} \ I get as
$\mu \le h^{-\nu}|\log h|^{-K}$ estimate (\ref{6-47}):

\begin{theorem}\label{thm-6-17} Under condition $(\ref{6-5})$ estimate
\begin{multline}
|\int \biggl( e(x,x, 0) -
(2\pi h)^{-1} \int {\rm e}_0(x_1,x_1,0;x_2,\xi_2)\,d\xi_2 \\-\cE^\MW (x,0) + \cE^\MW_0 (x,0)\biggr)
\psi (x)\,dx -
\sum \kappa_{l,m} \mu h^{1+2m} | \le C\mu^{-{1/\nu}}h^{-1}
\label{6-47}
\end{multline}
holds as $\mu \le \epsilon h^{-\nu}$.
\end{theorem}

\begin{proof} To finish the proof of this theorem one needs to cover the case
$h^{-\nu}|\log h|^{-K}\le \mu \le \epsilon h^{-\nu}$, getting rid of the term $|\log h|^K$ in the error estimates.

The first problematic error comes from the correction terms in proposition \ref{prop-6-15}, namely from the terms of the type
$\mu h^{1+2m}\int \varkappa_{l,m}(x_2)\gamma^{\nu -2m-3+k}\,dx$ with $k\ge 1$,
$\nu -2m-3+k=-1$ and this error term is $O(1)$ unless $k=1$, $\nu=2m+1$ in which case it it is $\kappa'_l \mu h^\nu |\log h|$. This is possible only for odd $\nu$ in which case operator $\cA_0$ is even with respect to $x_1\mapsto -x_1$,
$\xi_2\mapsto -\xi_2$ but perturbation contains exactly one factor $x_1$ and therefore it is odd and after integration with respect to $x_1$, $\xi_2$ this correction term results in 0 if $\psi$ is even with respect to $x_1$.

Further, one needs to consider terms corresponding to $m=3$ in the successive approximations leading to proposition \ref{prop-6-16} and there one can replace
$\cA _0(x_2,\xi_2)-\cA _0(y_2,\xi_2)$ by $B_1 (x_2-y_2)$, and also terms corresponding to $m=2$ in the same successive approximations and there one can replace $\cA _0(x_2,\xi_2)-\cA _0(y_2,\xi_2)$ by $B_2 (x_2-y_2)^2$.

To calculate the contribution of such terms one can apply the same approach as in the proof of proposition \ref{prop-6-15} and the contribution of $\gamma$-admissible partition element with respect to $x_1$ will be
\begin{equation*}
\sum _m \mu h^{1+2m} \int \varkappa_{l,m,k}(x_2) \psi (x) \gamma^{\nu-2m-3+k}\,dx
\end{equation*}
with $k\ge 0$; however since this expression should be $O(|\log h|^K)$ all the terms but those with $\nu\le 2m+1$, $k\ge 1$ should vanish; further, the total contribution of all remaining terms save those with $\nu = 2m+1$ and $k=1$ is $O(1)$, which leaves us with no ``bad'' terms for even $\nu$ and with one ``bad'' term $\kappa'_l \mu h^\nu \log h$ for odd $\nu$, $m=(\nu-1)/2$. However, parity considerations with respect to $x_1$ show that this term should vanish as well.
\end{proof}

\begin{remark}\label{rem-6-18-m}(i)
All the coefficients $\varkappa_{l,*}$ and $\kappa_{l,*}$  vanish for $l=0$.

\smallskip
\noindent
(ii) Obviously the same approximate expressions  (\ref{IRO6-3-52}), \ref{IRO6-3-52-*}, \ref{IRO6-3-52-**} as in \cite{IRO6} hold for part $\cE^\MW_\corr$ ``associated'' with $\cY_\inn$;
\end{remark}

\sect{Modified $V$. II. $\epsilon_0 \lowercase {h}^{-\nu}\le \mu \le C_0\lowercase { h}^{-\nu}$}

Now I will consider the intermediate case
\begin{equation}
\epsilon_0 h^{-\nu}\le \mu \le C_0 h^{-\nu}
\label{7-1}
\end{equation}
with arbitrarily small constant $\epsilon_0$ and arbitrarily large constant $C_0$; this case which described the largest possible values in \cite{IRO6}
now is no more than transition to the next section.

\subsection{Estimates}

Let us denote by
$\lambda_n(\xi_2)$ eigenvalues of operator
\begin{equation}
{\bf a}^0= {\frac 1 2}\biggl(D_1^2 + \bigl(\xi_2 - x_1^\nu/\nu\bigr)^2 - (2l+1)x_1^{\nu-1}\biggr);
\label{7-2}
\end{equation}
then $\Lambda_n(x_2,\xi_2)=\lambda_n (\xi_2) - {\frac 1 2} W(x_2)$ are eigenvalues of
${\bf a}= {\bf a}^0-W(x_2)$.

My main nondegeneracy assumption will be
\begin{equation}
|\Lambda_n|+ (|\xi_2|+1) |\partial_{\xi_2}\Lambda_n| +|\partial_{x_2}\Lambda_n|\ge \epsilon_0\qquad \forall n,\xi_2,
\label{7-3}
\end{equation}
may be coupled with \ref{6-4-pm}. This condition (\ref{7-3}) follows from
(\ref{6-5}); further, it follows from \ref{6-4-pm} for $|\xi_2|\ge C$.
On the other hand, since $\lambda_n\to 0$ and
$\xi_2\partial_{\xi_2}\lambda_n\to 0$ as $|\xi_2|\to \infty$, condition (\ref{7-3}) implies that $|W|+|\partial_{x_2}W|\ge \epsilon_0$ and therefore
locally one of conditions \ref{6-4-pm}, (\ref{6-5}) must be fulfilled. 

Obviously, under conditions (\ref{7-1}),(\ref{7-3}) for each $\xi_2$ number of eigenvalues of one-dimensional operator
\begin{equation}
\cA_0= {\frac 1 2}\biggl(h^2D_1^2 + \bigl(\xi_2 - \mu x_1^\nu/\nu\bigr)^2 - (2l+1)x_1^{\nu-1}-W\biggr)
\label{7-4}
\end{equation}
below level $c_0$ does not exceed $C$.

Further, note that condition (\ref{7-3}) for eigenvalues of $\cA_0$ is equivalent to the same condition for eigenvalues of ${\bf a}$. Then I easily arrive to

\begin{proposition}\label{prop-7-1} Under conditions $(\ref{7-1}),(\ref{7-3})$
contribution to the remainder estimate of the zone $\{|\xi_2|\le C\}$ is $O(1)$.
\end{proposition}

Furthermore, analysis in the zone $\cY''_0$ under condition (\ref{7-1}) does not differ from the analysis as $\mu\le \epsilon h^{-\nu}$. Namely

\begin{claim}\label{7-5}
Under conditions (\ref{7-1}) and \ref{6-4-pm} operator $\cA_0$ and thus operator $\cA$ is elliptic in the zone $\cY''_0\Def \{|\xi_2|\ge C\}$ and the contribution of $\cY''_0$ to the remainder estimate is negligible.
\end{claim}

\begin{claim}\label{7-6}
Similarly, under conditions (\ref{7-1}) and (\ref{6-5}) operator $\cA$ is microhyperbolic in the zone $\cY''_0\Def \{|\xi_2|\ge C\}$ and the contribution of $\cY''_0$ to the remainder estimate is $O(1)$.
\end{claim}

 Therefore

\begin{proposition}\label{prop-7-2} Let conditions $(\ref{7-1}),(\ref{7-3})$
and one of conditions \ref{6-4-pm}, $(\ref{6-5})$ be fulfilled. Then the remainder estimate is $O(1)$ where the principal part is defined by $(\ref{6-40})$.\end{proposition}

\subsection{Calculations} 
Calculations in this case also do not differ from those in section 6 leading to
the following statements

\begin{theorem}\label{thm-7-3}
Let conditions $(\ref{7-1})$, $(\ref{7-3})$ and \ref{6-4-pm} be fulfilled. Then $\cR_I$ defined by $(\ref{6-25})$ and $\cR^*$ defined by $(\ref{6-27})$ do not exceed $C$.
\end{theorem}

\begin{theorem}\label{thm-7-4}
Let conditions $(\ref{7-1})$ and $(\ref{6-5})$ be fulfilled. Then left-hand expressions of $(\ref{6-38})$, $(\ref{6-46})$ and $(\ref{6-47})$ do not exceed $C$.
\end{theorem}

\sect{Modified $V$. III. $ \mu \ge C_0\lowercase { h}^{-\nu}$}

Now I consider the previously forbidden case
\begin{equation}
\mu \ge C_0 h^{-\nu}
\label{8-1}
\end{equation}
with sufficiently large constant $C_0$. In this case all zones should be redefined. Also the difference between $l=0$ and $l\ge 1$ becomes crucial.

\subsection{Estimates. I}

As $|\xi_2|\asymp \mu \gamma^\nu$, $\gamma \ge C_1(\mu^{-1}h )^{1/(\nu+1)}$ let us consider first eigenvalues $\Lambda_n(x_2,\xi_2)$ of operator $\cA (x_2,\xi_2)$.
Then proposition \ref{prop-A-3} implies instantly that

\begin{claim}\label{8-2} As $n\ne l$ and $|\xi_2|\asymp \mu \gamma^\nu$, 
$\gamma \ge C_1(\mu^{-1}h )^{1/(\nu+1)}$
\begin{equation*}
\Lambda_n (x_2,\xi_2)  \asymp (n-l) \mu h \gamma ^{\nu-1}
\end{equation*}and signs of the left and right-hand expressions coincide and \begin{equation}\label{8-3}
\Lambda_l(x_2,\xi_2)=\omega_l h^2\gamma^{-2}-{\frac 1 2}W(x_2)+
O\Bigl(h^2\gamma^{-1} + h^2 (\mu^{-1}h)^2\gamma^{-4-2\nu}\Bigr), \quad
\omega_l > 0\ \text{as\;}\ l\ge1.
\end{equation}
\end{claim}\vskip-20pt
Therefore

\begin{claim}\label{8-4} As $l\ge 1$ zone
$\cY'''\Def \bigl\{  C_0(\mu h^\nu)^{1/(\nu+1)}\le |\xi_2|\le 
 \epsilon \mu h^\nu \bigr\}$ 
is elliptic and its contribution to the remainder estimate is $O(h^s)$.
\end{claim}
On the other hand,

\begin{claim}\label{8-5}
Under condition \ref{6-4-pm} zone $\cY''\Def \bigl\{|\xi_2|\ge C\mu h^\nu\bigr\}$ is elliptic as well and its contribution to the remainder estimate is  $O(h^s)$ as well for $l\ge 0$.
\end{claim}

Therefore as $l\ge 1$ and condition \ref{6-4-pm} is fulfilled, one needs to analyze only two remaining zones
$\cX_1=\bigl\{\epsilon \rho_1\le |\xi_2|\le C\rho_1 \bigr\}$, 
$\rho_1=\mu h^\nu$ and
$\cX_0=\bigl\{|\xi_2|\le C_0\rho_0 \bigr\}$, $\rho_0=(\mu h^\nu)^{1/(\nu+1)}$.

In the zone $\cX_1$ propagation speed with respect to $x_2$ is in average $\asymp \rho^{-1}$ (with $\rho=\rho_1$)  due to proposition \ref{prop-A-3} again and the propagation speed with respect to $\xi_2$ is in average $O(1)$ and therefore one can take
\begin{equation}
T_0=Ch|\log h|,\qquad T_1=\epsilon_1 \rho_1
\label{8-6}
\end{equation}
and for $T\in [T_0,T_1]$
propagation on the energy levels $\tau\in [-\epsilon,\epsilon]$ which started in $B(0,{\frac 1 2})$ does not leave $B(0,1)$ but the shift with respect to $x_2$
is $\asymp \rho^{-1}T$ and it satisfies logarithmic uncertainty principle and thus the spectral trace is negligible.

\begin{remark}\label{rem-8-1}
One should be more careful as $\mu \ge h^{-M}$ with arbitrarily large $M$ and use $\log \mu$ instead of $|\log h|$.
\end{remark}

Therefore
\begin{equation}
|F_{t\to h^{-1}\tau}{\bar \chi}_T(t) (Qu)|
\label{8-7}
\end{equation}
does not exceed $Ch^{-1}\rho T_0= C \rho |\log h|$ where $Q$ is a partition element corresponding to $\cX_1$, $|\tau|\le \epsilon$. Therefore due to Tauberian arguments  the contribution of this zone to the remainder is
$O(h^{-1}T_0 /T_1)=O(|\log h|)$. One can get rid off this superficial logarithmic factor both in the estimate of (\ref{8-7}) and in the remainder estimate; standard details I leave to the reader. So,

\begin{proposition}\label{prop-8-2} Let $l\ge 1$ and conditions \ref{6-4-pm} and $(\ref{8-1})$ be fulfilled. Then as $Q$ is supported in the zone $\cX_1$ expression $(\ref{8-7})$ does not exceed $C\rho_1$ and the contribution of $\cX_1$
to the remainder estimate is $O(1)$.
\end{proposition}

Therefore I am left with the zone
$\cX_0=\bigl\{|\xi_2|\le C_0(\mu h^\nu)^{1/(\nu+1)} \bigr\}$. Let us fix $x_2$. I don't know if eigenvalue $\lambda_n(\xi_2)$ of ${\bf a}^0(\xi_2)$ vanishes in $\cX_0$ (may be even with some of its derivatives)\footnote{\label{foot-10} It clearly happens for even $\nu$ and $n<l$.}
but I know that if it happens then $n\le c_1$; moreover due to the analyticity of $\lambda _n(\xi_2)$ it can happen only in no more then $C_1$ points and due to proposition \ref{prop-A-3}  and the analyticity of $\lambda _n(\xi_2)$  
\begin{equation}
\lambda _n (\eta)\sim \alpha (\eta-{\bar\eta})^r
\label{8-8}
\end{equation}
for some $\alpha \ne 0$ and $r=1,2,\dots$   near each such point ${\bar\eta}$, $\alpha$ and $r$ depend on ${\bar\eta}={\bar\eta}_{n,k}$ $k=1,\dots, K$ (depending on $\nu,l$ as well).  Further,   two eigenvalues do not vanish simultaneously.

But then condition \ref{6-4-pm} will provide non-degeneracy. Really, in our assumptions an ellipticity is broken only in the strips of the type
\begin{equation}
\cY=\bigl\{ |\xi_2 -{\bar\eta}\rho_0|\asymp C\Delta\bigr\},\qquad\Delta = \rho_0^{1-2/r},
\label{8-9}
\end{equation}
and the average propagation speed with respect to $x_2$ is of magnitude
$\rho_0^{-1}|\xi_2 -{\bar\eta}|^{r-1}\asymp \rho_0^{(2-r)/r}$ there
and therefore one can take
\begin{equation}
T_1=\epsilon \rho_0^{-(2-r)/r},\qquad
T_0 = Ch|\log h| \rho_0^{-(2-r)/r} \Delta^{-1}\asymp h|\log h|,\qquad \Delta = \rho_0^{1-2/r}.
\label{8-10}
\end{equation}
Therefore for $Q$ supported in the strip $\cY$ expression (\ref{8-7}) does not exceed $Ch^{-1}\Delta \times T_0 = C|\log h| \rho_0^{-(2-r)/r}$
and contribution of $\cY$ to the remainder estimate does not exceed this expression multiplied by $T_1^{-1}$ i.e. $Ch|\log h|$. Furthermore, using standard methods one can easily get rid off the superficial logarithmic factor both in the estimate of (\ref{8-7}) and the remainder estimate:

\begin{proposition}\label{prop-8-3} Let $l\ge 1$ and conditions \ref{6-4-pm} and $(\ref{8-1})$ be fulfilled. Then as $Q$ is supported in the strip $\cY$ described by $(\ref{8-9})$, expression $(\ref{8-7})$ does not exceed $C\rho_0^{-(2-r)/r}$ and the contribution of $\cY$ to the remainder estimate is $O(1)$. 
\end{proposition}

Therefore I arrive to

\begin{proposition}\label{prop-8-4}
Let $l\ge 1$ and conditions \ref{6-4-pm} and $(\ref{8-1})$ be fulfilled. Then the remainder estimate is $O(1)$ while the principal part is given by $(\ref{6-40})$ for different strips with any $T\in [T_0,T_1]$ defined by $(\ref{8-10})$ for strip $\cY$ under conditions $(\ref{8-8})-(\ref{8-9})$ and by $(\ref{8-6})$ for strip $\cX_1$.
\end{proposition}

I would like to note that

\begin{proposition}\label{prop-8-5} Let $l\ge 1$ and conditions $(\ref{6-4})_-$ and $(\ref{8-1})$ be fulfilled. Then 

\smallskip
\noindent
(i) Zone $\cX_1$ is elliptic and its contribution to the remainder estimate is $O(h^s)$;

\smallskip
\noindent
(ii) Furthermore if also condition
\begin{equation}
\lambda_n(\eta)\ne 0 \qquad \forall n,\eta
\label{8-11}
\end{equation}
is fulfilled\,\footnote{\label{foot-11} However I cannot check condition (\ref{8-11}).} then the remainder estimate is $O(h^s)$.
\end{proposition}

\subsection{Estimates. II}

Let us consider the special case $l=0$; I remind that  then only eigenvalue $\lambda_0(\eta)$ should be considered and that condition $(\ref{6-4})_-$ leads then to the asymptotics with the principal part 0 and remainder estimate $O(h^s)$ and therefore is excluded from the further consideration.

Further, as $\nu$ is odd $\lambda_0=0$ identically, condition $(\ref{6-4})_+$ provides ellipticity everywhere. Thus I arrive to

\begin{proposition}\label{prop-8-6}
Let $l=0$, $\nu$ be odd and conditions $(\ref{6-4})_+$ and $(\ref{8-1})$ be fulfilled. Then the remainder estimate is $O(h^s)$ while the principal part is given by $(\ref{6-40})$.
\end{proposition}

On the other hand, as $l=0$, $\nu$ is even and condition $(\ref{6-4})_+$ holds due to proposition \ref{prop-A-7} ellipticity is violated only in the strip
\begin{equation}
\cY=\{ \epsilon_1 \Delta \le |\xi_2 - \eta \rho_0 | \le C\Delta\}, \qquad \eta \asymp |\log \rho_0|^{\nu/(\nu+1)},\;
\Delta =\rho_0 |\log \rho_0|^{-1/(\nu+1)}
\label{8-12}
\end{equation}
where as before $\rho_0 =(\mu h^\nu)^{1/(\nu+1)}$. In this strip propagation speed with respect to $x_2$ is $\asymp \Delta^{-1}$ and again
\begin{equation}
T_0= Ch|\log h|,\qquad T_1= \epsilon \Delta
\label{8-13}
\end{equation}
and expression (\ref{8-7}) does not exceed
$Ch^{-1}\Delta T_0 = C\Delta |\log h|$ and the remainder estimate is
$O(|\log h|)$. Further, by the standard arguments one can get rid off the superficial logarithmic factors. Thus

\begin{proposition}\label{prop-8-7}
Let $l=0$, $\nu$ be even and conditions $(\ref{6-4})_+$ and $(\ref{8-1})$ be fulfilled. Then the remainder estimate is $O(1)$ while the principal part is given by $(\ref{6-40})$ with $T_0,T_1$ defined by $(\ref{8-13})$.
\end{proposition}

\subsection{Estimates. III}

Now I want to derive estimates under condition \ref{6-4-pm} replaced by (\ref{6-5}). Without condition \ref{6-4-pm} some zones cease to be elliptic and should be reexamined:

\begin{claim}\label{8-14}
As $l\ge 1$ these zones are
$\bigl\{|\xi_2|\ge C \mu h^\nu\bigr\}$ and also\end{claim}\vskip-20pt
\begin{claim}\label{8-15}
As $l\ge 1$ these zones are ``inner parts'' of the strips  described by (\ref{8-9}), namely,
$\cY=\bigl\{ |\xi_2 - {\bar\eta }\rho_0 | \le \epsilon_1\Delta\bigr\}$.
\end{claim}\vskip-20pt
\begin{claim}\label{8-16} As $l=0$, $\nu $ even this zone is
$\bigl\{|\xi_2|\ge C\rho_0 |\log \rho_0|^{\nu/(\nu+1)}\bigr\}$;
\end{claim}\vskip-20pt
\begin{claim}\label{8-17} As $l=0$, $\nu$ odd this zone is 
$\bigl\{|\xi_2|\le \epsilon \mu\bigr\}$.
\end{claim}

Since condition (\ref{6-5}) provides $T_0=Ch|\log h|$ anyway contribution of (\ref{8-9})-type strips to the remainder estimate will be $O(1)$ again. The standard partition-rescaling arguments in all other zones bring contribution of all other zones to $O(\log \mu)$; however additional arguments of the proof of proposition \ref{prop-6-14} allow us to reduce it to $O(1)$. Therefore

\begin{proposition}\label{prop-8-8} Let conditions $(\ref{8-1})$ and $(\ref{6-5})$ be fulfilled. Then
the remainder estimate is $O(1)$ while the principal part of the asymptotics is given by $(\ref{6-40})$ for different zones with any $T\in [T_0,T_1]$, $T_0=Ch|\log h|$ and $T_1$ defined as in propositions \ref{prop-8-2}--\ref{prop-8-7}.
\end{proposition}

\subsection{Calculations. I}
In this subsection I give the principal parts of asymptotics already derived under condition \ref{6-4-pm} in more explicit form.

First of all, consider method of successive approximations fixing $x_2=y_2$. Then while contribution of the strip of the width
$\Delta $ in $\xi_2$ to the principal part is of magnitude $\Delta h^{-1}$, each next term of successive approximations acquires factor
$|\partial_{\xi_2}\Lambda _n| T \times T/h\asymp
(\partial_{\xi_2}\Lambda _n) h|\log h|^2$ with $T=T_0$ where $\Lambda_n$ is an eigenvalue of $\cA$.
Further one needs to consider only strips where ellipticity fails and then $\Delta \asymp |\partial_{\xi_2}\Lambda_n|^{-1}$.

So, the first, the second and the the third terms do not exceed
\begin{pdeq}\label{8-18}\end{pdeq}
\begin{equation}
Ch^{-1}|\partial_{\xi_2}\Lambda _n|^{-1},\quad
C |\log h|^2, \quad
C h |\partial_{\xi_2}\Lambda_n |\cdot |\log h|^4
\tag*{$(\theequation)_{1-3}$}\label{8-18-*}
\end{equation}
respectively.

Actually the second term in the successive approximations is $O(1)$. Really, considering the second term which corresponds to the linear part $(x_2-y_2)\partial_{y_2}\cA(y_2,hD_2)$ of the perturbation one can rewrite it 
as the result of direct calculations in the form including 
$\partial_{x_2}\partial_{\xi_2}\Lambda _n =0$; on the other hand considering the second term corresponding to the rest $(x_2-y)^2\cB (x_2,y_2,hD_2)$ of the perturbation one can estimate it easily by $O(h^\delta)$.

Now I can rewrite the principal part of the asymptotics as
\begin{equation}
(2\pi h)^{-1}\int {\rm e}(x_1,x_1,0;x_2,\xi_2)\psi (x)\,d\xi_2 dx
\label{8-19}
\end{equation}
with error not exceeding already achieved remainder estimate which is either $O(1)$ or $O(h^\infty)$ (where remainder estimate $O(h^\infty)$ corresponds to the elliptic case and no successive approximations are needed at all).

Let us consider the contribution of the strips where ellipticity is broken to the error; I remind it does not exceed the minimum of all three expressions in \ref{8-18-*}. Then $(\ref{8-18})_3$ is obviously $O(1)$ in all cases with the singular exception of the strip  (\ref{8-9}) with $r=1$, 
$\rho h \ge |\log h|^{-K}$. However in this case $(\ref{8-18})_1$ is $O(1)$
unless $|\log h|^{-K}\le \rho h \le |\log h|^K$ and one can still handle this case getting rid off the superficial logarithmic factors in $(\ref{8-18})_{1,3}$  by the standard arguments. Thus I arrive to

\begin{theorem}\label{thm-8-9} Let conditions \ref{6-4-pm} and $(\ref{8-6})$ be fulfilled. Then

\smallskip
\noindent
(i) Asymptotics with the principal part given by $(\ref{8-19})$ holds with the remainder estimate $O(1)$;

\smallskip
\noindent
(ii) Furthermore, as $l=0$, $\nu$ is odd this asymptotics holds with the remainder estimate $O(h^\infty)$.
\end{theorem}

Furthermore, fixing $W$ at $x_1=0$ and $\alpha=1 $ and thus replacing $\cA$ by $\cA^0$ to the pilot model operator, I can apply the method of successive approximation again; then each next term  gets an extra factor
$C\gamma T_0 h^{-1} |\log h|$ with $\gamma=(\mu ^{-1}|\xi_2|)^{1/\nu}$ and only strips where ellipticity breaks should be counted. Also one can see easily that

\begin{claim}\label{8-20}
The error does not exceed the second term
$Ch^{-2}T_0 \Delta \gamma$\,\footnote{\label{foot-12} I skip superficial logarithmic factors one can easily get rid off by the standard arguments.}.
Furthermore, for odd $\nu$ and perturbation, which is odd with respect to $x_1$, the second term is 0 and therefore the error does not exceed the sum of the second term with a perturbation $O(x_1^2)$ and the third term with a perturbation $O(x_1)$ i.e. $Ch^{-3}T_0^2 \Delta \gamma^2$\,$^{\ref{foot-12}}$.
\end{claim}

Thus, I just list the different cases:

\begin{claim}\label{8-21}
As $l\ge 1$ and condition $(\ref{6-4})_+$ is fulfilled the main contribution to the error is provided by the zone $\cX_1$ with $\xi_2 \asymp \mu h^\nu$ and
$\gamma \asymp h$ and of the width
$\Delta \asymp \mu h^\nu$; so the error is $O\bigl(\mu h^\nu\bigr)$. The contributions of (\ref{8-9})-type strips  are much smaller;
\end{claim}
\vskip-10pt
\begin{claim} \label{8-22} As $l\ge 1$ and condition $(\ref{6-4})_-$ is fulfilled the main contribution to the error is provided by (\ref{8-9})-type strips with the largest possible $r$; then $\xi_2 =O\bigl( (\mu h^\nu)^{1/(\nu+1)}\bigr)$,
$\gamma \asymp (\mu ^{-1}h)^{1/(\nu+1)}$ and
$\Delta \asymp (\mu h^\nu)^{(r-2)/r(\nu+1)}$; so the error is
$O\bigl( (\mu h^\nu)^{-\delta}\bigr)$ with $\delta = 2/r(\nu+1)$ anyway;
\end{claim}
\vskip-10pt
\begin{claim} \label{8-23} As $l= 0$, $\nu$ is even and condition $(\ref{6-4})_+$ is fulfilled the main contribution to the error is provided by  $\cX_1$ with
$\xi_2 \asymp (\mu h^\nu)^{1/(\nu+1)}|\log (\mu h^\nu)|^{\nu/(\nu+1)}$,
$\gamma \asymp (\mu ^{-1}h)^{1/(\nu +1)} |\log (\mu h^\nu)|^{1/(\nu+1)}$ and of the width
$\Delta \asymp (\mu h^\nu )^{1/(\nu+1)}|\log (\mu h^\nu)|^{-1/(\nu+1)}$;
so the error is $O(1)$ anyway;
\end{claim}
\vskip-10pt
\begin{claim} \label{8-24} As $l= 0$, $\nu$ is odd and condition $(\ref{6-4})_+$ is fulfilled the error is just $O(h^\infty)$.
\end{claim}

Thus I arrive to asymptotics with the principal part
\begin{equation}\label{8-25}
(2\pi h)^{-1}\int {\rm e}_0(x_1,x_1,0;x_2,\xi_2)\psi (x)\,d\xi_2 dx
\end{equation}
and remainder estimates described in Theorem \ref{thm-8-10} below:

\begin{theorem}\label{thm-8-10} Let condition $(\ref{8-1})$ be fulfilled. Then

\smallskip
\noindent
(i) As $l\ge 1$ and condition $(\ref{6-4})_+$ is fulfilled asymptotics with the principal part given by $(\ref{8-25})$ holds with the remainder estimate 
$O(\mu h^\nu)$;

\smallskip
\noindent
(ii) As \underline{either} $l\ge 1$ and condition $(\ref{6-4})_-$ is fulfilled \underline{or} $l=0$, $\nu$ is even and condition $(\ref{6-4})_+$ is fulfilled asymptotics with the principal part given by $(\ref{8-25})$ holds with the remainder estimate $O(1)$;

\smallskip
\noindent
(iii) Furthermore, as $l=0$, $\nu$ is odd and condition $(\ref{6-4})_+$ is fulfilled the same asymptotics holds with the remainder estimate $O(h^\infty)$.
\end{theorem}

\subsection{Calculations. II}

In this subsection I give in more explicit form the principal parts of asymptotics already derived under condition (\ref{6-5}). Basically I need
to reconsider only the external formerly elliptic zones described by 
(\ref{8-14})--(\ref{8-17}). The analysis in the first of them is not different from the analysis under condition \ref{6-4-pm}; analysis in the second one repeats the proof of theorem \ref{thm-6-17}; analysis in two latter is rather obvious. Thus I arrive to two following theorems:

\begin{theorem}\label{thm-8-11} Let conditions $(\ref{6-5})$ and $(\ref{8-6})$ be fulfilled. Then asymptotics with the principal part $(\ref{8-19})$
holds with the remainder estimate $O(1)$.
\end{theorem}

\begin{theorem}\label{thm-8-12} Let conditions $(\ref{6-5})$ and $(\ref{8-6})$ be fulfilled. Then

\smallskip
\noindent
(i) As $l\ge 1$ estimate
\begin{multline}
\cR^{**}\Def|\int \biggl( e(x,x, 0) -
(2\pi h)^{-1} \int {\rm e}_0(x_1,x_1,0;x_2,\xi_2)\,d\xi_2 \\-\cE^\MW (x,0) + \cE^\MW_0 (x,0)\biggr)
\psi (x)\,dx -
\sum \kappa_{l,m} \mu h^{1+2m} | \le C\mu h^\nu
\label{8-26}
\end{multline}
holds;

\smallskip
\noindent
(ii) As $l=0$ estimate
\begin{multline}
|\int \biggl( e(x,x, 0) -
(2\pi h)^{-1} \int {\rm e}_0(x_1,x_1,0;x_2,\xi_2)\,d\xi_2 \\-\cE^\MW (x,0) + \cE^\MW_0 (x,0)\biggr)
\psi (x)\,dx  | \le C
\label{8-27}
\end{multline}
holds.
\end{theorem}

\appendix\sect{Appendix: Eigenvalues of 1D operators}

\subsection{General observations}

In this Appendix $\lambda_n(\eta)$ ($n=0,1,\dots$) denote eigenvalues of one-dimensional pilot-model Schr\"odinger operators with $\mu=h=1$
\begin{align}
&{\bf a}^0(\eta) = D ^2 +(\eta - x^\nu/\nu)^2 - (2l+1) x^{\nu-1}
\label{A-1}\\
\intertext{or more general operator}
&{\bf a} (\eta) =(1+\alpha_1 x+ \beta_1^2x^2) D ^2 +
(1+\alpha_2 x+ \beta_2^2x^2)(\eta - x^\nu/\nu)^2 -\label{A-2} \\
&\qquad\qquad\qquad\qquad\qquad\qquad\qquad\qquad\qquad\qquad\qquad\qquad(2l+1) (1+\alpha_3 x)x^{\nu-1}
\notag
\end{align}
with $\nu=2,3,\dots$ and $\beta_j> \alpha_j^2/2$.

One can prove easily the following   statement:

\begin{proposition}\label{prop-A-1} Let $l\in \bR$. Then

\smallskip
\noindent
(i) As $|\eta |\le C_0$ the spacing between two consecutive eigenvalues $\lambda_n$ and $\lambda_{n+1}$ with $n\le c_0$ is
$\asymp 1$;

\smallskip
\noindent
(ii) For  operator $(\ref{A-1})$ with odd $\nu$ \ $\lambda_n(-\eta)=\lambda_n(\eta)$;

\smallskip
\noindent
(iii) For even $\nu$ and $\eta \le 0$ \
$\lambda_n (\eta)\ge (1-\epsilon)\eta^2 - C_1$ \; $\forall n=0,1,\dots$.
\end{proposition}

However, the case of even $\nu$ and $\eta\to -\infty$ is rather exceptional:

\begin{proposition}\label{prop-A-2} As
$\eta \ge C_0$ (and thus also as $\eta \le -C_0$ and $\nu$ is odd)

\smallskip
\noindent
(i) The spacing between eigenvalues with $n\le c_0$ is
$\asymp (1+|\eta|)^{(\nu-1)/\nu}$;

\smallskip
\noindent
(ii) As $n< l$ \ ($l< n \le c_0$)
$\lambda_n(\eta)$ is less than (greater than respectively)
\newline $\epsilon (n-l) (1+|\eta|)^{(\nu-1)/\nu}$\,\footnote{\label{foot-?}
Thus leaving the special case $n=l\in \bZ^+$ for the further analysis.}.
\end{proposition}

\begin{proof} Proof follows from the proof of proposition \ref{prop-A-3} below.

\end{proof}

\subsection{Asymptotic behavior of $\lambda_l(\eta)$ as
$\eta\to \infty$ as $l\ge1$}

In this subsection I prove

\begin{proposition}\label{prop-A-3} (i) For operator $(\ref{A-1})$ with $l\ge 1$ as $\eta \to +\infty$
(and thus also as $\eta \to -\infty$ and $\nu$ is odd)
\begin{equation}
\lambda_l(\eta)= \kappa \eta^{-2/\nu}+ O\bigl( \eta^{-(\nu+3)/\nu}\bigr)
\label{A-3}
\end{equation}
with $\kappa > 0$;

\smallskip
\noindent
(ii) For operator $(\ref{A-2})$ with $l\ge 1$ as $\eta \to +\infty$
(and thus as $\eta \to -\infty$ and $\nu$ is odd)
\begin{equation}
\partial_{\alpha_j} \lambda_l(\eta)\bigr|_{\alpha = \beta =0} =
\kappa_j \eta + O(\eta^{-1/\nu})
\label{A-4}
\end{equation}
with $\kappa_1=\kappa_2=-\kappa_3/2$, $\alpha=(\alpha_1,\alpha_2,\alpha_3)$,
$\beta=(\beta_1,\beta_2,\beta_3)$ and furthermore
\begin{equation}
\sum_{1\le j\le 3}\partial_{\alpha_j} \lambda_l(\eta)\bigr|_{\alpha = \beta =0} =
\kappa_4 \eta ^{1/\nu}\lambda_l + O(\eta^{-2/\nu}).
\label{A-5}
\end{equation}
\end{proposition}

\begin{proof} (i) Let us plug $\eta =\gamma^\nu/\nu$ with $\gamma\gg 1$ where in the case even $\nu$ this is the only scenario and in the case of odd $\nu$ analysis of scenario $\xi_2=-\gamma^\nu/\nu$ is done by the symmetry. Then after shift
$x\mapsto x+\gamma$ operator ${\bf a}^0(\eta)$ is transformed into operator
\begin{multline*}
D^2 + x^2\Bigl( \gamma^{\nu-1}+ {\frac 1 2}(\nu-1)x \gamma^{\nu-2}+
{\frac 1 6}(\nu-1)(\nu-2)x^2\gamma^{\nu-3}+ \dots \Bigr)^2\\
- (2l+1)\Bigl(\gamma^{\nu-1} + (\nu-1)x \gamma^{\nu-2}+
{\frac 1 2}(\nu-1)(\nu-2) x^2\gamma^{\nu -3} +\dots\Bigr)\
\end{multline*}
and after rescaling $x\mapsto x\gamma^{(1-\nu)/2}$ this operator is transformed into $\gamma^{\nu-1}{\bf b}_\varepsilon$ where
\begin{multline*}
{\bf b}_\varepsilon= D^2 + x^2\Bigl( 1+
{\frac 1 2}(\nu-1)x \varepsilon +
{\frac 1 6}(\nu-1)(\nu-2)x^2\varepsilon^2 + \dots \Bigr)^2\\
- (2l+1)\Bigl(1 + (\nu-1)x \varepsilon+
{\frac 1 2}(\nu-1)(\nu-2) x^2\varepsilon^2 +\dots \Bigr)
\end{multline*}
with $\varepsilon=\gamma^{-(\nu+1)/2}$. Then
\begin{multline*}
{\bf b}_\varepsilon= \underbrace{D^2 + x^2 -(2l+1)}_{{\bf h}_0}
+\varepsilon \underbrace{(\nu-1) \Bigl( x^3 -(2l+1)x\Bigr)}_{{\bf h}_1} +\\
\varepsilon^2\underbrace{(\nu-1)\Bigl( ({\frac 7 {12}}\nu-{\frac {11}{12}}) x^4 -{\frac 1 2}(2l+1)(\nu-2)x^2\Bigr)}_{{\bf h}_2} +O\bigl(\varepsilon^3\bigr)
\end{multline*}
and let us denote by $\Lambda_\varepsilon$ and $U_\varepsilon$ its eigenvalue close to $0$ and the corresponding eigenfunction. Then
\begin{equation}
\Lambda _\varepsilon = \omega_1\varepsilon +\omega_2 \varepsilon^2 +\dots\qquad\text{and}\qquad
U_\varepsilon = u _0 + u_1\varepsilon + u_2\varepsilon^2 \dots
\label{A-6}
\end{equation}
where obviously $u_0=\upsilon_l$ is a Hermite function, $\omega_1=\omega_3=\dots=0$
and
\begin{equation}
{\bf h}_0 u_1 + {\bf h}_1 u_0=0\qquad
{\bf h}_0 u_2 +{\bf h}_1u +{\bf h}_2u_0= \omega_2 u_0.
\label{A-7}
\end{equation}
Then
\begin{equation}
\omega_2 = \langle {\bf h}_1u +{\bf h}_2u_0, u_0\rangle=
-\langle u , {\bf h}_0u\rangle +\langle {\bf h}_2u_0, u_0\rangle.
\label{A-8}
\end{equation}
It is known that $(x-iD)\upsilon_k=(2k+2)^{1/2}\upsilon_{k+1}$,
$(x+iD)\upsilon_k=(2k)^{1/2}\upsilon_{k-1}$ and therefore
\begin{align*}
&x \upsilon_l = &&{\frac 1 2}\Bigl(
(2l+2)^{1/2}\upsilon_{l+1} + (2l)^{1/2}\upsilon_{l-1}\Bigr),\\
&x^2 \upsilon_l = &&{\frac 1 4}\Bigl(
(2l+2)^{1/2}(2l+4)^{1/2}\upsilon_{l+2} + (4l+2)\upsilon_l +
(2l)^{1/2}(2l-2)^{1/2}\upsilon_{l-2}\Bigr),\\
&\bigl(x ^2-2l-1)\upsilon_l = &&{\frac 1 4}\Bigl(
(2l+2)^{1/2}(2l+4)^{1/2}\upsilon_{l+2} -
2(2l+1) \upsilon_l+
(2l)^{1/2}(2l-2)^{1/2}\upsilon_{l-2}\Bigr),\\
&x\bigl(x ^2-2l-1)\upsilon_l = &&{\frac 1 8}\Bigl(
(2l+2)^{1/2}(2l+4)^{1/2}(2l+6)^{1/2}\upsilon_{l+3} -
(2l+2)^{1/2}(2l-2)\upsilon_{l+1}-\\
&&&\hphantom{ {\frac 1 2}\Bigl( }(2l)^{1/2}(2l+4)\upsilon_{l-1}+
(2l)^{1/2}(2l-2)^{1/2}(2l-4)^{1/2}\upsilon_{l-3}\Bigr),
\end{align*}
which imply
\begin{multline*}
\hskip-\labelsep\langle {\bf h}_0 u,u\rangle = \\{\frac 1 {64}}(\nu-1)^2\Bigl( {\frac 1 6}(2l+2)(2l+4)(2l+6) + {\frac 1 2} (2l+2)(2l-2)^2 -{\frac 1 2} (2l)(2l+4)^2 -{\frac 1 6} (2l) (2l-2)(2l-4)\Bigr)=\\
{\frac 1 {16}} (\nu-1)^2 \Bigl(-2l^2-2l+3\Bigr).
\end{multline*}
On the other hand
\begin{multline*}
\langle {\bf h}_2u_0,u_0\rangle = (\nu-1)\Bigl(
({\frac 7 {12}}\nu-{\frac {11}{12}}) \|x^2u_0\|^2 -
{\frac 1 2}(\nu-2)(2l+1)\|xu_0\|^2\Bigr)= \\
(\nu-1)({\frac 7 {12}}\nu-{\frac {11}{12}}){\cdot}{\frac 1 {16}}\Bigl(
(2l+2)(2l+4)+(4l+2)^2+(2l)(2l-2)\Bigr) -\\
{\frac 1 4}(\nu-1)(\nu-2)\cdot (2l+1)^2=\\
(\nu-1)(7\nu-11){\cdot}{\frac 1 {16}}\Bigl(2l^2+2l+1\Bigr) -
{\frac 1 4}(\nu-1)(\nu-2)(2l+1)^2
\end{multline*}
and
\begin{multline*}
\omega_2= {\frac 1 {16}} (\nu-1)\biggl( (7\nu-11)\bigl(2l^2+2l+1\bigr) -
4(\nu-2)\bigl(4l^2+4l+1\bigr)- (\nu-1) \bigl(-2l^2-2l+3\bigr)\biggr)=\\
{\frac 1 2} (\nu-1)l (l+1),
\end{multline*}
Therefore $\Lambda _\varepsilon = \omega_2 \varepsilon^2+O(\varepsilon^4)$ as $\varepsilon \to 0$ (because $\omega_3=0$ as well) which implies statement (i) with $\kappa =\omega_2\nu^{-2/\nu}$.

\smallskip
\noindent
(ii) After obvious transformations
\begin{equation*}
\partial_{\alpha_j} \lambda_l(\eta)\bigr|_{\alpha = \beta =0} =
\gamma^{\nu-1} \langle {\bf k}_j U_\varepsilon, U_\varepsilon\rangle
\end{equation*}
with
\begin{align*}
&{\bf k}_1 = (\gamma + \varepsilon x )D^2,\\
&{\bf k}_2= x^2\Bigl( 1+
{\frac 1 2}(\nu-1)x \varepsilon +
{\frac 1 6}(\nu-1)(\nu-2)x^2\varepsilon^2 + \dots \Bigr)^2\\
&{\bf k}_3=- (2l+1)\Bigl(1 + (\nu-1)x \varepsilon+
{\frac 1 2}(\nu-1)(\nu-2) x^2\varepsilon^2 +\dots \Bigr)
\end{align*}
and therefore
\begin{equation*}
\langle {\bf k}_j U_\varepsilon, U_\varepsilon\rangle = \gamma
\langle {\bf k}'_j u_0,u_0\rangle +O\bigl(\varepsilon^2\gamma\bigr)
\end{equation*}
with ${\bf k}'_1=D^2$, ${\bf k}'_2=x^2$, ${\bf k}'_3=-(2l+1)$
which implies (\ref{A-4}). 

Known equalities
$\langle x^2\upsilon_l,\upsilon_l\rangle =
\langle D^2\upsilon_l,\upsilon_l\rangle = (2l+1)/2$ imply that
$\kappa_1=\kappa_2=-\kappa_3/2$. Further,
$\sum_{1\le j\le 3} \langle {\bf k}_jU_\varepsilon, U_\varepsilon\rangle
= \gamma \lambda_l +O(\varepsilon^2)$ which implies (\ref{A-5}).
\end{proof}

\subsection{More general operators}

Now I consider operator
\begin{equation}
\cA (y,\eta)\Def \beta\Bigl(\alpha h^2D^2 \alpha + \alpha ^{-2}(\eta -\mu x^\nu/\nu)^2 -(2l+1)\mu h x^{\nu-1}\Bigr)\beta
\label{A-9}
\end{equation}
with
\begin{equation}
\alpha =\alpha (x,y),\quad \beta=\beta(x,y) ,\quad \alpha (0,y)=1,\quad
c_0^{-1}\le \beta \le c_0.
\label{A-10}
\end{equation}
Let $\lambda_n$ be eigenvalues of $\cA$. Changing 
$x\mapsto \gamma (\mu^{-1}h)^{1/(\nu+1)}x$ and
$\eta \mapsto (\mu h^\nu)^{1/(\nu+1)}$ respectively  I arrive to operator (\ref{A-9})  again  with
$\mu = h=1$ and $\alpha$, $\beta$ replaced by
$\alpha ((\mu^{-1}h)^{1/(\nu+1)}x, y)$, $\beta ((\mu^{-1}h)^{1/(\nu+1)}x, y)$ and with a factor  $(\mu h^\nu)^{2/(\nu+1)}$.

\begin{proposition}\label{prop-A-4} Let conditions $(\ref{A-9}),(\ref{A-10})$ be fulfilled. Then

\smallskip
\noindent
(i) $\lambda_n(\eta)\ge C_0(\mu h^\nu)^{2/(\nu+1)}$ as $n\ge C$;

\smallskip
\noindent
(ii) As $|\eta |\le C_0(\mu h^\nu)^{2/(\nu+1)}$ the spacing between consecutive eigenvalues with $n\le c_0$ is $\asymp (\mu h^\nu)^{2/(\nu+1)}$ and
\begin{equation}
|\partial_y^p \partial_\eta ^q\lambda_n(y,\eta)|\le C_{pq} (\mu^{-1}h)^{p/(\nu+1)}
(\mu h^\nu)^{(2-q)/(\nu+1)};
\label{A-11}
\end{equation}

\smallskip
\noindent
(iii) For even $\nu$ and $\eta \le 0$ \quad
$\lambda_n(y,\eta)\ge (1-\epsilon) \eta^2 -C_1$, $n=0,1,\dots$
\end{proposition}

\begin{proposition}\label{prop-A-5}
As $\eta \ge C_0(\mu h^\nu)^{1/(\nu+1)}$ (and thus also as
$\eta \le -C_0(\mu h^\nu)^{1/(\nu+1)}$ and  $\nu$ is odd)

\smallskip
\noindent
(i) The spacing between eigenvalues with $n\le c_0$ is
$\asymp |\eta|^{(\nu-1)/\nu} (\mu h^\nu)^{1/\nu}$;

\smallskip
\noindent
(ii) As $n<l$ \ ($l<n \le c_0$) $\lambda_n(y,\eta)$ is less than (greater than respectively) \newline
$\epsilon (n-l)\bigl((\mu h^\nu)^{2/(\nu+1)} + |\eta|^{(\nu-1)/\nu} (\mu h^\nu)^{1/\nu}\bigr)$ 
and these eigenvalues satisfy
\begin{equation}
|\partial_y^p \partial_\eta ^q\lambda_n(y,\eta)|\le C_{pq}(\mu^{-1}h)^{p/(\nu+1)}
|\eta |^{-q}|\lambda_n(y,\eta)|;
\label{A-12}
\end{equation}

\smallskip
\noindent
(iii) As $\eta \ge C_0(\mu h^\nu)$ (and thus as
$\eta \le -C_0(\mu h^\nu)$ and  $\nu$ is odd) 
$|\lambda_l(y,\eta)|\le \epsilon_0$.
\end{proposition}

An extra analysis is needed for our purposes as $n=l$ and
\begin{equation}
\mu h^\nu \ge C_1
\label{A-13}
\end{equation}
with arbitrarily large $C_1$.

\begin{proposition}\label{prop-A-6}
Let condition $(\ref{A-13})$ be fulfilled and $l\ge 1$. Then
as $\eta \ge C_0(\mu h^\nu)^{1/(\nu+1)}$ 
\begin{equation}
\lambda_l (y,\eta)\asymp (\mu h^\nu /\eta)^{2/\nu}\qquad\text{and}\qquad
\eta \partial_\eta \lambda_l (y,\eta)\asymp (\mu h^\nu /\eta)^{2/\nu}.
\label{A-14}
\end{equation}
\end{proposition}

\subsection{Case of $\lambda_l$ as $l=0$}

Here cases of  odd and even $\nu$ differ drastically. Note first that
\begin{equation}
{\bf a}^0 (\eta)= \bigl(i D + \xi_2 - x^\nu/\nu\bigr)
\bigl(-i D + \xi_2 - x^\nu/\nu\bigr)
\label{A-15}
\end{equation}
and as $\nu$ is odd operator ${\bf a}^0(\eta)$ has the bottom eigenvalue $\lambda_0(\eta)$ with eigenfunction defined from
$\bigl(-\partial + \xi_2 - x^\nu/\nu\bigr)v=0$ i.e. $v=\exp (\xi_2 x - x^{\nu+1}/\nu(\nu+1))$ and therefore $\lambda_0(\eta)$ is identically $0$.

Similarly, as $\beta=1$ operator $\cA$ defined by (\ref{A-9}) is equal modulo $O(h^2)$ to operator
\begin{multline}
\cB(y,\eta)\Def h^2 \alpha ^2 D +\alpha^{-2} (\eta -\mu x^\nu/\nu)^2- \mu h x^{\nu-1}=\\
\bigl(ihD \alpha +\alpha^{-1} (\eta -\mu x^\nu/\nu)\bigr)
\bigl(-\alpha ihD +\alpha^{-1} (\eta -\mu x^\nu/\nu)\bigr)
\label{A-16}
\end{multline}
and I arrive to the statement (i) of

\begin{proposition}\label{prop-A-7} (i) For odd $\nu$ the bottom eigenvalue of $\cB(y,\eta)$ is $0$;

\smallskip
\noindent
(ii) For even $\nu$ the bottom eigenvalue of $\cB(y,\eta)$ is
$(\mu h^\nu)^{2/(\nu+1)} \Lambda \bigl(y,\eta (\mu h^\nu)^{-1/(\nu+1)}\bigr)$
where
\begin{align}
&C^{-1}\exp (-C\eta^{(\nu+1)/\nu}) \le \Lambda (y,\eta) \le
C \exp (-\epsilon \eta^{(\nu+1)/\nu}),\label{A-17}\\
&\epsilon \eta ^{1/\nu}\le -\partial_\eta \bigl(\log \Lambda (y,\eta)\bigr) \le C \eta^{1/\nu}.
\label{A-18}
\end{align}
\end{proposition}

\begin{proof} I need to consider the case of even $\nu$ only. The same representation (\ref{A-15}) shows that $\lambda_0 (y,\eta)>0$. However, since this eigenfunction is fast decaying outside of the potential well, one can do
the same shift and rescaling as before and using arguments of \cite{HeMa} to prove that
$\Lambda_0 (y,\eta) \sim k\exp (-k_2 \eta^{(\nu+1)/\nu})$. Also one can prove easily
that $\partial_\eta \Lambda_0 (y,\eta)\sim -k_3 \eta ^{1/\nu} k\exp (-k_2 \eta^{(\nu+1)/\nu})$ as $\eta \ge C$ with $k_3=kk_2(1+\nu)/\nu$. Estimates 
(\ref{A-17}), (\ref{A-18}) follow from this.
\end{proof}

\input IRO7.bbl

\end{document}

%% file: IRO7.bbl
\bibliographystyle{amsalpha}

\providecommand{\bysame}{\leavevmode\hbox to3em{\hrulefill}\thinspace}

\vglue .06truein

\hfill\hfill {\sl   \today \/}

\vglue .06truein

\begin{tabular}{rrl}
&{\hskip 220 pt} &Department of Mathematics,\cr
&&University of Toronto,\cr
&&100, St.George Str.,\cr
&&Toronto, Ontario M5S 3G3\cr
&&Canada\cr
&&ivrii@math.toronto.edu\cr
&&Fax: (416)978-4107\cr
\end{tabular}